\documentclass[lettersize,journal]{IEEEtran}
\usepackage{cite}
\usepackage{color}
\usepackage{amsfonts,amssymb}
\usepackage{bm}
\usepackage{amsmath}
\usepackage{multirow}
\usepackage{algorithmic}
\usepackage{graphicx}
\usepackage{subfigure}
\usepackage{float}
\usepackage{textcomp}
\usepackage{graphicx} 
\usepackage{epstopdf}
\usepackage{array}
\usepackage[thmmarks,amsmath]{ntheorem}
\newtheorem{mythm}{Theorem}
\newtheorem{mydef}{Definition}
\newtheorem{assumption}{Assumption}
\newtheorem{lemma}{Lemma}
\newtheorem{remark}{Remark}
\newtheorem{proposition}{Proposition}

\newtheorem{corollary}{Corollary}
\newtheorem{problem}{Problem}

\usepackage{array}
\usepackage{enumerate}
\usepackage{booktabs}
\usepackage{algorithm}
\usepackage{algorithmic}
\usepackage{setspace}
\usepackage{float}
\usepackage{cases}

\usepackage{mathrsfs}

\begin{document}

\title{Global solution to sensor network localization: 
	\\ A non-convex potential game approach and  its distributed implementation }

\author{Gehui Xu, Guanpu Chen,  Yiguang Hong, \IEEEmembership{Fellow, IEEE}, Baris Fidan, \IEEEmembership{Senior Member, IEEE},  
\\
Thomas Parisini, \IEEEmembership{Fellow, IEEE}, Karl H. Johansson, ~\IEEEmembership{Fellow,~IEEE}
\thanks{
This work was supported by the Swedish Research Council, the Knut and Alice Wallenberg Foundation, and the Swedish Foundation for Strategic Research, and was also supported 
in part by the Digital Futures Scholar-in-Residence Program, by the European Union's Horizon 2020 research and innovation programme under grant agreement no. 739551 (KIOS CoE), and by the Italian Ministry for Research in the framework of the 2017 Program for Research Projects of National Interest (PRIN), Grant no. 2017YKXYXJ.
}
	\thanks{Gehui Xu is with Key Laboratory of Systems and Control, Academy of Mathematics and Systems Science, Beijing 100190, China,
(e-mail: xghapple@amss.ac.cn)}
\thanks{Guanpu Chen and Karl H. Johansson  are with the Division of Decision and Control Systems, School of
	Electrical Engineering and Computer Science, KTH Royal Institute of
	Technology,  and also with Digital Futures, SE-10044
Stockholm, Sweden (e-mail: guanpu@kth.se,  kallej@kth.se)}
 	\thanks{Yiguang Hong is with
		Department of Control Science and Engineering, Tongji University,  Shanghai 201804, China, and is also with Shanghai Research Institute for Intelligent Autonomous Systems, Shanghai 201210 China. (e-mail: yghong@iss.ac.cn)}
 \thanks{Baris Fidan is with the Department of Mechanical and Mechatronics Engineering, University of Waterloo, Waterloo, ON N2L 3G1, Canada (e-mail:
fidan@uwaterloo.ca)}
   	\thanks{ Thomas Parisini is with  the Department of Electrical and Electronic Engineering,
	Imperial College London, London SW7 2AZ, UK, and also with the Department
	of Engineering and Architecture, University of Trieste, Trieste 34127, Italy. (e-mail: t.parisini@imperial.ac.uk)}
}



\maketitle
{
\begin{abstract}
Consider a sensor network consisting of both anchor and non-anchor nodes. We address the following sensor network localization (SNL) problem:
given the physical locations of anchor nodes and relative measurements among all nodes, determine the locations of all non-anchor nodes.
The solution to the  SNL problem is challenging due to its inherent non-convexity. 
In this paper, the problem takes on the form of a multi-player non-convex potential game in which canonical duality theory is used to define a complementary dual potential function.
After showing the Nash equilibrium (NE) correspondent to the  SNL solution, we provide a  necessary and sufficient condition for a stationary point to coincide with the NE. An algorithm is proposed to reach the NE and shown to have convergence rate  $\mathcal{O}(1/\sqrt{k})$.
With the aim of reducing the information exchange within a network, a distributed algorithm for NE seeking is implemented and its global convergence analysis is provided.
Extensive simulations show the validity and effectiveness of the proposed approach to solve the SNL problem.
\end{abstract}}


\section{Introduction}

Wireless sensor networks (WSNs) -- owing to their sensing, processing, and communication capabilities -- have a wide range of applications \cite{akyildiz2002wireless,liu2015event}, for instance, target tracking and detection \cite{marechal2010joint,jing2021angle}, environment monitoring \cite{sun2010corrosion}, area exploration \cite{sun2005reliable},  and cooperative robot tasks \cite{jing2018weak}.
		For all these applications, it is essential to accurately determine the location of every sensor. 
Estimating locations of the sensor nodes based on measurements between neighboring nodes has attracted significant research interest in recent years (see, for example, \cite{estrin2000embedding,savvides2001dynamic}). 

Range-based methods constitute a common inter-node measurement approach utilizing transmission-based techniques such as
time of arrival, time-difference of arrival, and strength of received radio frequency signals  \cite{mao2007wireless}.
Due to limited transmission power, the measurements can only be obtained within a radio range. A pair of nodes are called neighbors if their distance is less than this radio range \cite{wan2019sensor}. 
Also, in many cases, there are  anchor nodes whose global positions are known
\cite{ahmadi2021distributed}. 

The sensor network localization
(SNL) problem is defined as follows:  Given the positions of the WSN  anchor nodes and the relative measurements among each non-anchor node and its neighbors,  find the positions of all nodes.

\subsection{Motivations and Challenges}

To  describe individual preferences and network interactions, game theory constitutes an effective modeling tool 
\cite{jia2013distributed,bejar2010cooperative}. The Nash equilibrium (NE) is a key solution concept characterizing a profile of stable strategies in which rational players would not choose to deviate from their own strategies~\cite{nash1951non,xu2022efficient}. The SNL problem can be formulated as a game by letting non-anchor nodes be players, their estimated position be strategies, and their positioning accuracy measurements be payoffs.
Specifically,  potential games are well-suited to model 
SNL  (see, for instance, \cite{jia2013distributed,moragrega2015potential,ke2017distributed}).
In particular, the potential game paradigm guarantees an alignment between the individual sensor profits and the network objective by exploiting a  unified potential function.   It is possible to find the  NE corresponding to the global optimum for the whole WSN rather than a local approximation.

{It is worth noting, though, that non-convexity is an intrinsic challenge of SNL problems, which, unfortunately, cannot be avoided by potential games or other modeling approaches. 
Indeed, finding the global solution to the SNL problem is still open, both in centralized \cite{wang2006further,tseng2007second} and distributed settings \cite{wan2019sensor,deghat2011distributed,calafiore2010distributed}.

In a centralized setting, some relaxation methods such as semidefinite programming (SDP) \cite{wang2006further}  and second-order cone programming (SOCP) \cite{tseng2007second} are employed to transform the original   non-convex problem into a convex one,  and then determine the positions of all unknown sensor nodes through a data fusion center. These relaxation methods ignore the non-convex  constraints, yielding only an approximate solution. 
    Some distributed methods have been developed to find a global network solution using only local
information,  such as alternating rank minimization (ARMA)-based algorithms \cite{wan2019sensor}, Barycentric     coordinate-based algorithms \cite{deghat2011distributed},  and gradient-based algorithms with Barzilai-Borwein stepsizes  \cite{calafiore2010distributed}.   These studies also fail to adequately address the non-convexity,   either lacking the theoretical guarantee of convergence or relying on additional prior initialization.
   }
 \vspace{-0.2cm}
\subsection{Contributions}  
{The aim of this paper is to seek solutions to the SNL problem and develop the corresponding algorithms. We  
 formulate the non-convex SNL problem as a potential game and 
 investigate the algorithm's global convergence to an NE via a conjugate transformation.
 Moreover, to address the communication capacity limitations and complexity  issues, we propose a distributed implementation of the NE-seeking algorithm.
 }
The contributions of this paper are summarized as follows:
\begin{enumerate}[(i)]
{	\item 
	We formulate the non-convex SNL problem as a multi-player potential game, with the individual payoff function and the overall potential function defined in terms of fourth-order polynomial forms, in order to effectively deal with the non-convex structures.
Then, we show that the NE corresponds to the global SNL solution. 

	\item We construct a conjugate transformation of the non-convex SNL potential game via canonical duality theory. 
Then, we provide a  necessary and sufficient condition on gradients of the components of the potential function for reaching the NE.  With this guidance, 
we propose a conjugate-based algorithm to compute the  NE  and establish that the convergence rate of this algorithm is  $\mathcal{O}(1/\sqrt{k})$.

 %
	
	
	\item 
To reduce the transmission  of dual information in the network,
we design a distributed implementation of the conjugate-based algorithm, borrowing inspiration from sliding mode control and extra-gradient techniques.
We show that the distributed algorithm can globally converge to the NE.
} 
\end{enumerate}

{A preliminary short version of this manuscript was submitted to 2024 European Control Conference 
\cite{Xu2023nonconvex}.
Different from the results in \cite{Xu2023nonconvex}, which only deals with a sufficient condition to reach the  NE and a centralized conjugate-based algorithm,   the current paper studies the necessary condition for reaching the NE,  the equivalency between the NE and the global SNL solution, and the convergence rate of the conjugate-based algorithm. Furthermore, the current paper provides a distributed SNL algorithm, together with its convergence analysis.
}

\vspace{-0.2cm}
\subsection{Related work}
\vspace{-0.15cm}
    In this section, we briefly review SNL game-based approaches in the literature, as well as methods as to how to handle the non-convexity issue and the main distributed SNL algorithms.
    

\noindent\textbf{Game-based approaches to  SNL } 
 {
Different game models have been developed to represent SNL problems, including 
 coalition games \cite{bejar2010cooperative}, zero-sum games \cite{vempaty2012localization},  and potential games \cite{jia2013distributed}. 
 Potential games are used in describing SNL problems in \cite{jia2013distributed,moragrega2015potential,ke2017distributed}.
 Specifically, {the so-called consensus of cognition among anchors and the consensus of physical measurements were considered  for the localization problem in \cite{ke2017distributed},  utilizing the time-of-arrival technique to characterize the potential game.}    
 A distributed potential game was proposed in  \cite{jia2013distributed} for the localization of road sensor networks,  employing received signal strength (RSS) based distance measurement and a quadratic distance estimation error cost.
A similar approach was  adopted in \cite{moragrega2015potential}    to construct the potential function for a distributed power control and node selection problem,  employing an RSS measurement together with a lognormal path loss model. Despite these efficient formulations leading to fruitful achievements, they only obtain an approximate NE or a local NE. Thus, these approaches do not guarantee reaching the global solution of the SNL problem.
 }

\noindent\textbf{Handling non-convexity in SNL} 
Range-based SNL is usually formulated as a nonlinear equality constrained feasibility problem.
 Relaxation methods, such as SDP \cite{wang2006further,lui2008semi}  and SOCP \cite{srirangarajan2008distributed}, were applied to this formulation by relaxing the original problem into a computationally tractable convex optimization problem.  These relaxation methods usually ignore the non-convex rank constraints,  leading only to an approximate solution to the SNL problem. In this view,  ARMA  was considered to obtain an exact solution by equivalently converting the rank constraints into complementary constraints \cite{wan2019sensor}, only guaranteeing local convergence.
 Another formulation of the range-based SNL is to solve a non-convex least squares problem. In this setting, \cite{calafiore2010distributed}  considered the gradient-based method to find the optimal solution while \cite{jia2013distributed} employed a stochastic best-response scheme to obtain the optimal solution of each iteration. However, these works do not deal with the non-convex structures, so most of their proposed algorithms lead to local optima dependent on the initial conditions.
 

%
%
%

\noindent\textbf{Distributed implementation of  SNL} 
Distributed algorithms that localize all sensor nodes with local information exchanged among neighboring nodes have gained attention in recent years.
SDP-based  \cite{wang2006further,biswas2006distributed}  and  SOCP-based \cite{srirangarajan2008distributed,shi2010distributed} are two classical methods extended  from centralized to distributed approaches for the SNL problem.  For example,
 \cite{biswas2006distributed}  divided the overall sensor nodes into clusters and solved the SDP problem at each cluster.
   \cite{shi2010distributed} established the communication network for exchanging the estimation information among sensor nodes and proposed the  computation using  an edge-based SOCP in each node. 
However, these distributed algorithms are based on    relaxation approaches,  again only yielding an approximation solution of the problem.     An ARMA-based distributed algorithm was proposed in \cite{wan2019sensor},
where the original SNL problem was decomposed into a group of node-based ARMA subproblems. Here, each subproblem had only one node to be localized and all other sensors were fixed, using the results obtained from the previous iteration. Still, this method lacks a theoretical guarantee of global convergence.
 Barycentric coordinate-based distributed algorithms \cite{deghat2011distributed} and  gradient-based distributed  algorithms  with Barzilai-Borwein stepsizes  \cite{calafiore2010distributed} have been developed to find the global optimal solution of the SNL problem with additional assumptions.  Nevertheless, \cite{deghat2011distributed}  required that the sensors should be located in the convex hull of their neighbors, and \cite{calafiore2010distributed} required that the sensors should be initialized with sufficiently good prior knowledge of the node positions. 

\vspace{-0.25cm}
\subsection{Outline}
\vspace{-0.15cm}
 The paper is organized as follows: Section II introduces the preliminaries. Section  III formulates a non-convex SNL potential game, while Section IV  investigates the  derivation of the NE  followed by the construction of an SNL algorithm.
 Section~\ref{sec:dis}  develops 
a distributed implementation for seeking  NE and gives its  convergence analysis.
Following that, 
 Section~\ref{sim} examines the effectiveness of the
 proposed approach with several experiments.
  Finally, Section VII concludes the paper
  and offers some future directions.

\section{Preliminaries}\label{pre}

\noindent\textbf{Notation} 
Let $ \mathbb{R}^{n} $(or $ \mathbb{R}^{m\times n} $) denote the set of $ n $-dimensional (or $ m $-by-$ n $) real column vectors (or real matrices), 
$ \boldsymbol{1}_{n} $(or $ \boldsymbol{0}_{n} $) denote the $ n $-dimensional column vector with all elements of $ {1} $ (or $ {0} $), 
 $ \operatorname{col}\{x_{1}\!,\dots\!,x_{n}\} \!  \!=\! \!(\!x^{\mathrm{T}}_{1},\! \dots,\!x^{\mathrm{T}}_{n}\!)^{\mathrm{T}} $, $\|\cdot\|$ denote the Euclidean norm of vectors,  and  $ \nabla f $ denote  function $  f $'s gradient. 

\noindent\textbf{Convex Analysis}
A set $ K \subseteq \mathbb{R}^{n} $ is convex if $ \omega x_{1}+(1-\omega)x_{2} \in K$ for  any
$ x_{1}, x_{2}\in K $ and $ 0\leq\omega\leq 1   $. 
For a closed convex set $ K $, the projection map $ \Pi_{K}: \mathbb{R}^{n} \rightarrow K $ is defined as
$
\Pi_{K}(x) \triangleq {\operatorname{argmin}}_{y \in K}\|x-y\|.
$
The following two basic results hold \cite{facchinei2003finite}:
\begin{equation}\label{projection}
	\begin{aligned}
		& \left(x-\Pi_{K}(x)\right)^\mathrm{T}\left(\Pi_{K}(x)-y\right) \geq 0, \quad \forall y \in K, \\
		& \left\|\Pi_{K}(x)-\Pi_{K}(y)\right\| \leq\|x-y\|, \quad \forall x, y \in \mathbb{R}^n.
	\end{aligned}
\end{equation}
A mapping $ F : \mathbb{R}^{n} \rightarrow \mathbb{R}^{n} $ is said to be  monotone  on a set $ K\subseteq \mathbb{R}^{n} $ if 
$
(F(x)-F(y))^{\mathrm{T}}(x-y) \geq 0, \quad \forall x, y \in K.
$

\noindent\textbf{Duality Theory}
  A differentiable function $\Psi: \Theta\rightarrow \mathbb{R}$ is said to be a {canonical function} on $\Theta  \subseteq \mathbb{R}^{q}$ if  its
derivative
  $\nabla\Psi : \Theta\rightarrow \Theta^{*}\subseteq \mathbb{R}^{q}$ is  one-to-one and onto from $\Theta$ to its range $\Theta^{*}$. 	
%
Besides, if $\Psi$ is a convex canonical function, then its 
conjugate function $\Psi^{*}: \Theta^{*} \rightarrow \mathbb{R} $ is uniquely defined by the Legendre transformation as $\Psi^{*}\left(\sigma\right)=\left\{\xi^{\operatorname{T}} \sigma-\Psi(\xi) \mid \sigma=\nabla \Psi(\xi), \xi\in \Theta\right\}$,
where 	$\sigma\in \Theta^{*} $
is a canonical dual variable.
Moreover, 	the following canonical duality relation holds on $\Theta\times\Theta^{*}$ \cite{gao2017canonical}: 
\begin{align*}
	\sigma&=\nabla \Psi(\xi)\Leftrightarrow\xi=\nabla \Psi^{*}\left(\sigma\right)\Leftrightarrow \xi^{\operatorname{T}} \sigma= \Psi(\xi)+\Psi^{*}\left(\sigma\right)
 \vspace{-0.3cm}
\end{align*} 	
$(\xi,\sigma)$ is called the Legendre canonical duality pair
on $\Theta\times\Theta^{*}$.

{In our approach, we consider both graph and game representations of the sensor network and SNL of interest. 
Here, denote $ \mathcal{N}= \mathcal{N}_{s}\cup \mathcal{N}_{a}$ as the sensor node set, where  $ \mathcal{N}_{s} $ and $ \mathcal{N}_{a}$ represent the  non-anchor and anchor node set, respectively.   In the game representation, the elements of $\mathcal{N}_{s}$ are also regarded as players.}

{
\noindent\textbf{Graph Theory}
{An undirected graph is denoted  by $\mathcal G(\mathcal N, \mathcal E)$ where $\mathcal N=\{i\}_{i=1}^{\bar{N}}$ corresponds to the node (or vertex) set  and $\mathcal E\subseteq\mathcal N \times \mathcal N$ is the edge set.  
{A $n$-dimensional ($n\in\{2,3\}$) representation of graph is a mapping of  $\mathcal G(\mathcal N, \mathcal E)$  to the point 
formations $\bar{\boldsymbol{x}}: \mathcal N \rightarrow \mathbb{ R}^{n}$, where 
$\bar{\boldsymbol{x}}(i)=x_i^{\operatorname{T}}$ is  the row vector of the coordinates of
the $i$-th node in $\mathbb{ R}^{n}$ and $x_i\in \mathbb{ R}^{n}$.  In the  SNL problem, the  $x_i$ is the position of sensor node $i$, and 
$ \mathcal{E}=\{(i,j)\in\mathcal{N}\times\mathcal{N}:\|x_{i}-x_{j}\|\leq R_s,i\neq j\} \cup  \{(i,j)\in\mathcal{N}_a\times\mathcal{N}_a:i\neq j\}$ is the edge set,  i.e., there is an edge between two sensor nodes if and only if either they are neighbors based on a transmission range $R_s$ or they are both anchors, where $R_s$ represents a sensor’s capability of sensing range measurements from others.
Given a graph  $\mathcal G(\mathcal N, \mathcal E)$ and an $n$-dimensional representation $\bar{\boldsymbol{x}}$ of it, the pair $(\mathcal G, \bar{\boldsymbol{x}})$ is called a $n$-dimensional framework. A framework  $(\mathcal G, \bar{\boldsymbol{x}})$ is called generic\footnote{Some special configurations exist among the sensor positions, e.g., groups of sensors may be collinear. The reason for using the term generic is to highlight the need to exclude the problems arising from such configurations.} if the set containing the coordinates of all its points is algebraically independent over the
rationales \cite{anderson2010formal}.  A framework $(\mathcal G, \bar{\boldsymbol{x}})$ is called  rigid if there exists a sufficiently small positive constant $\epsilon$ such that if every framework $(\mathcal G, \bar{\boldsymbol{y}})$
satisfies  $\|x_i-y_i \|\leq\epsilon$ for $i\in\mathcal N$ and 
$\| x_i-x_j\|=\|y_i-y_j\|$ for every  pair $i,j\in \mathcal N$ connected by an edge in $\mathcal E$, then 
$\| x_i-x_j\|=\|y_i-y_j\|$ holds  for any node pair $i,j\in \mathcal N$  no matter there is an edge between them. A   graph
$\mathcal G(\mathcal N, \mathcal E)$  is called generically   $n$-rigid or simply rigid (in $n$ dimensions) if any generic framework 
$(\mathcal G, \bar{\boldsymbol{x}})$ is rigid. A framework $(\mathcal G, \bar{\boldsymbol{x}})$  is globally rigid if every framework $(\mathcal G, \bar{\boldsymbol{y}})$ satisfying
$\| x_i-x_j\|=\|y_i-y_j\|$  for any node pair $i,j\in \mathcal N$  connected by an edge in $\mathcal E$ and  $\| x_i-x_j\|=\|y_i-y_j\|$  for any node pair $i,j\in \mathcal N$  that are not connected by a single edge. 
 A graph
$\mathcal G(\mathcal N, \mathcal E)$  is called generically globally rigid if 
any generic framework 
$(\mathcal G, \bar{\boldsymbol{x}})$  is 
globally rigid \cite{anderson2010formal,fidan2010closing,tay1985generating}.}}


\noindent\textbf{Game Theory} 
A  standard potential game-theoretical model can be denoted by $G=\{\mathcal{N}_{s}, \{\Omega_i\}_{i\in\mathcal{N}_{s}}, \{J_{i}\}_{i\in\mathcal{N}_{s}}\}$.
 For each player $i\in \mathcal{N}_{s}=\{1,\dots,N\}$,   the $i$-th player has a strategy variable $x_i$  in a local feasible set  $\Omega_i\in\mathbb{R}^n$.
 Denote $   \boldsymbol{\Omega}\triangleq\prod_{i=1}^{N}\Omega_{i} \subseteq \mathbb{R}^{nN} $, $ \boldsymbol{x}\triangleq \operatorname{col}\{x_{1}, \dots, x_{N}\} \in \boldsymbol{\Omega} $ as the strategy profile for all players,  and $ \boldsymbol{x}_{-i}\triangleq \operatorname{col}\{x_{1}, \dots ,x_{i-1}, x_{i+1}, \dots, x_{N}\}$ as the strategy profile for all players except player $ i $. Player $i$ has a  continuously differentiable payoff function
 $J_i(x_i, \boldsymbol{x}_{-i}): \mathbb{R}^{nN}\rightarrow \mathbb{R}$,
 and aims to minimize its own payoff by taking all players’ strategies into account. 
Each player needs to consider others' strategies when choosing its optimal strategy.  The best-known concept therein
is the NE \cite{nash1951non}.
\begin{mydef}[NE]\label{d1}
	A profile $ \boldsymbol{x}^{\Diamond}=\operatorname{col}\{x_{1}^{\Diamond}, \dots, x_{N}^{\Diamond} \} \in \boldsymbol{\Omega} \subseteq \mathbb{R}^{nN}$ is said to be an  NE of game (\ref{f1}) if 
	\begin{equation}\label{ne}
		J_{i}\left(x_{i}^{\Diamond}, \boldsymbol{x}_{-i}^{\Diamond}\right) \leq J_{i}\left(x_{i}, \boldsymbol{x}_{-i}^{\Diamond}\right),\;\forall i \in \mathcal{N}_s, \forall x_{i}\in \Omega_{i}. \vspace{-0.2cm}
	\end{equation}
\end{mydef}
The NE above characterizes a strategy profile $\boldsymbol{x}^{\Diamond}$ that each player adopts its optimal strategy.  
Given others' strategies $\boldsymbol{x}_{-i}^{\Diamond}$, no one can benefit from changing its strategy unilaterally.

Then we  introduce the concept of potential game
\cite{monderer1996potential}. 
\begin{mydef}[potential game]\label{d2}
	A game $G=\{\mathcal{N}_{s}, \{\Omega_i\}_{i\in\mathcal{N}_{s}}, \{J_{i}\}_{i\in\mathcal{N}_{s}}\}$ is a potential
	game if there exists a potential function $\Phi$ such that, for $i\in\mathcal{N}_s$,
	\begin{equation}\label{pp1}
	\Phi(x_{i}^{\prime}, \boldsymbol{x}_{-i})-\Phi\left(x_{i}, \boldsymbol{x}_{-i}\right)=J_{i}(x_{i}^{\prime}, \boldsymbol{x}_{-i})-J_{i}\left(x_{i}, \boldsymbol{x}_{-i}\right), \vspace{-0.2cm}
	\end{equation}
	for every  $\boldsymbol{x}\in \boldsymbol{\Omega}$,  and unilateral deviation $x_{i}^{\prime}\in \Omega_i$. 
\end{mydef}
It follows from Definition   \ref{d2} that 
any unilateral deviation from a strategy profile always results in the same change in both individual payoffs and a unified potential function.}



\section{Problem formulation}

In this section, we first introduce the range-based SNL problem and then formulate it as a potential game.

{
\subsection{Sensor network localization problem}
	 Consider a static sensor network 
	 in $\mathbb{R}^{n}$  ($ n = 2  $ or 3)  composed of 
	  $ M $ anchor nodes
whose positions are known  and $N$ non-anchor  nodes whose positions are
unknown (usually $ M<N $). 
The network of sensors can be indexed by $ \mathcal{N}\!=\! \mathcal{N}_{s}\cup \mathcal{N}_{a}$, where  $ \mathcal{N}_{s}\!=\!\{1,2,\!\dots\!,N\} $ and $ \mathcal{N}_{a}\!=\! \{N\!+\!1,N\!+\!2,\dots,N\!+\!M \} $ are the sets of non-anchor nodes and anchor nodes, respectively.  
{Let $ x_{i}^{\star} \in \mathbb{R}^{n} $ for $i\in \mathcal{N}$ denote the {\color{blue} } position of the $ i $-th  sensor, where $x_{j}$ for $j\in \mathcal{N}_s$ is the  non-anchor node's position, and	
$ a_{l}=x_{N+k}^{\star}$  
for $l\in \mathcal{N}_{a}$ and $k\in \{1,2,\dots,M \}$
is the anchor node's position.
}
For a pair of non-anchor node $ i\in \mathcal{N}_s$ and anchor node $ l\in \mathcal{N}_a $, their Euclidean distance
is denoted as $ e_{il} $.
Similarly, the distance between non-anchor node  $i $ and
non-anchor node  $ j $ is denoted as $ d_{ij} $. 
Consider a graph $ \mathcal{G}=(\mathcal{N}, \mathcal{E}) $  representing the sensing relationships between sensors. Here,  $\mathcal{E}=\mathcal{E}_{ss}\cup \mathcal{E}_{as}\cup \mathcal{E}_{aa}$ with 
$ \mathcal{E}_{ss}=\{(i,j)\in\mathcal{N}_{s}\times\mathcal{N}_{s}:\|x_{i}^{\star}-x_{j}^{\star}\|\leq R_s,i\neq j\} $  denoting the edge set between non-anchor nodes,  	$ \mathcal{E}_{as}=\{(i,l)\in\mathcal{N}_{s}\times\mathcal{N}_{a}:\|x_{i}^{\star}-a_{l}\|\leq R_s\} $ denoting the edge set between anchor nodes and non-anchor nodes, and $ \mathcal{E}_{aa}=\{(l,m)\in\mathcal{N}_{a}\times\mathcal{N}_{a},l\neq m\} $ denoting the edge set between anchor nodes.
Denote $\mathcal{N}_{s}^{i}\!=\!\{j:(i,j)\!\in\!\mathcal{E}_{ss}  \}$ as  the neighbor set of adjacent non-anchor nodes 
for non-anchor  node $i$,
and $\mathcal{N}_{a}^{i}\!=\!\{l:(i,l)\!\in\!\mathcal{E}_{as}  \}$ as  the  neighbor set of adjacent anchor nodes 
for  non-anchor node $i$.
Also, two sensor nodes can communicate with each other if and only if the distance between them is smaller than a communication range $R_c$. Here we consider $R_s=R_c$.

We make the following basic assumption.
\begin{assumption}
 The sensor topology graph	 $ \mathcal{G} $ is undirected and generically globally rigid.
\end{assumption}
The undirected graph topology is a widely-used assumption in many graph-based approaches \cite{chen2021distributed,jing2021angle}. 
The connectivity of $\mathcal{G} $ can also be induced by some disk graphs \cite{wan2019sensor},  which ensures the validity of the information transmission between nodes.
{The generic global rigidity of $\mathcal{G} $ has been widely employed in SNL problems to guarantee the graph structure invariant, which indicates a unique localization of sensor network   \cite{calafiore2010distributed,shi2010distributed,anderson2008rigid}.
Besides, there have been extensive discussions on graph rigidity in existing works \cite{wan2019sensor,cao2021bearing}, but it is not the primary focus of our paper. }

Let us now define the SNL  problem.
Suppose that 
the measurements $d_{ij}$ and $e_{il}$ are noise-free and all anchors’ positions $ a_{l} $, $l\in \mathcal{N}_{a}$ are accurate. 
 The
{ range-based SNL task is to determine the accurate positions  of all non-anchor sensor nodes, $i\in \mathcal{N}_{s}$ when all anchor node positions $ a_{l} $, $l\in \mathcal{N}_{a}$ and  measurements $d_{ij}$ and $e_{il}$ are given, 
that is,  }
\begin{align}\label{SNL}
	& \text{find} \quad x_{1},\dots,x_{N}\in \mathbb{R}^{n}\\
	& \text{s.t.} \quad \|x_{i}-x_{j}\|=d_{ij}, \forall (i,j)\in \mathcal{E}_{ss},\notag\\
	& \quad\quad\;\|x_{i}-a_{l}\|=e_{il}, \forall (i,l)\in \mathcal{E}_{as}. \notag \vspace{-0.2cm}
\end{align}
}
 {\color{blue}
 }
{Denote $\boldsymbol{x}^{\star}=\operatorname{col}\{x_{1}^{\star}, \dots, x_{N}^{\star} \} \in \mathbb{R}^{nN} $ as the actual position vector of all non-anchor nodes. Assumption 1 guarantees that (\ref{SNL}) has a unique solution and is equal to $\boldsymbol{x}^{\star}$ \cite{eren2004rigidity}. }
Then we 
 formulate
 a multi-player non-convex potential game problem to
 reach the solution $\boldsymbol{x}^{\star} $ of non-convex (\ref{SNL}) and develop an implementable algorithm
for this solution in the sequel.

\subsection{Potential game formulation}

In the  SNL problem, each non-anchor node needs to consider the location accuracy of the whole sensor network while ensuring its own positioning accuracy through the given information. 
In other words, each non-anchor node needs to ensure the alignment between its individual objective and collective objective, which will be defined in the sequel. To this end, we formulate the non-convex SNL (\ref{SNL})  as a potential game model to describe non-anchor nodes' individual preferences and their interactions within the network \cite{jia2013distributed,monderer1996potential}. Also, the  NE of the potential game is verified to be the network’s precise localization through a potential function. 
{Define $G=\{\mathcal{N}_{s}, \{\Omega_i\}_{i\in\mathcal{N}_{s}}, \{J_{i}\}_{i\in\mathcal{N}_{s}}\}$ as an $N$-player SNL potential game, where $\mathcal{N}_{s}=\{1,\dots,N\}$ corresponds to the player set,
$\Omega_i$ is player $i$'s local feasible set, which is convex and compact, and $J_{i}$ is player $i$'s payoff function. In this context, we map the position estimated by each non-anchor node as each player's chosen strategy. i.e.,  the strategy of the player $i$ (non-anchor node) is the estimated position $x_i\in \Omega_i$.  
For $i\in\mathcal{N}_s$, the payoff function $J_{i}$ is constructed as
 \begin{equation*}\label{jI}
 J_{i}(x_{i},\boldsymbol{x}_{-i})\!=\!\sum_{j\in\mathcal{N}_{s}^{i} }(\|x_{i}-x_{j}\|^{2}\!-d_{ij}^2)^2+\!\sum_{l\in\mathcal{N}_{a}^{i} }(\|x_{i}-a_{l}\|^{2}-e_{il}^2)^2,  \vspace{-0.2cm}
 \end{equation*}
  where the first term in $J_{i}$  measures the   localization accuracy between  non-anchor node $i$ and its  non-anchor node neighbor $j\in\mathcal{N}_{s}^{i}$ and  the second term  measures the  localization accuracy  between  $i$ and its   anchor neighbor $l\in\mathcal{N}_{a}^{i}$. 

 The  individual objective of each non-anchor node  is to ensure its position accuracy,
 i.e.,
\begin{equation}\label{f1}
	\min \limits_{x_{i} \in \Omega_{i}} J_{i}\left(x_{i}, \boldsymbol{x}_{-i}\right). \quad 
\end{equation}
%
To achieve each non-anchor node's accurate localization,  the concept of NE $ \boldsymbol{x}^{\Diamond}$  is considered according to Definition \ref{d1}. The NE reflects that each non-anchor node's estimated position is optimal when all non-anchor nodes attain  $ \boldsymbol{x}^{\Diamond}$. 
Here, we sometimes write an NE as a \textit{global} NE to 
distinguish our setting  from \textit{local} NE \cite{jia2013distributed,ke2017distributed,heusel2017gans}, which  only satisfies condition \eqref{ne} within a small neighborhood $\Omega_{i}\bigcap\delta(x_{i}^{\Diamond})$
for $i \in \mathcal{N}_s$, rather than the whole $\Omega_{i}$.
Moreover, 
by regarding the individual payoff  $J_{i}$ as a marginal contribution to the whole network's collective objective \cite{marden2009cooperative,jia2013distributed}, we consider 
the following measurement of 
 the overall performance of the sensor network:  
\begin{equation}\label{potential-fun}
	\Phi\!\left(
 \boldsymbol{x}
 \right)\!=\!\!\!\!\!\!\!\sum_{(i,j)\in\mathcal{E}_{ss} }\!\!\!\!\!(\|x_{i}-x_{j}\|^{2}\!-d_{ij}^2)^2\!+\!\!\!\!\!\!\sum_{(i,l)\in\mathcal{E}_{as} }\!\!\!\!\!(\|x_{i}-a_{l}\|^{2}-e_{il}^2)^2.
 \vspace{-0.15cm}
\end{equation}
    {Here, $J_{i}$ denotes the   localization accuracy of non-anchor node $i$, which depends on the  strategies of $i$'s neighbors, while  $	\Phi$ denotes the   localization accuracy of the entire network according to \eqref{SNL}. 

The following lemma reveals the relationship between the solutions of (\ref{potential-fun}) and \eqref{SNL} with the proof given in Appendix \ref{r122} 

\begin{lemma}\label{r11}
{ Under Assumption 1, the global minimum of $\Phi$ in (\ref{potential-fun}) is unique and is equal to  $\boldsymbol{x}^{\star}$ of \eqref{SNL}.}
\end{lemma}


Then we verify that  $\Phi$ is a   potential function in Definition \ref{d2}, thereby  ensuring the existence of global NE $\boldsymbol{x}^{\Diamond}$  and its consistency with $\boldsymbol{x}^{\star}$. The proof is shown in Appendix \ref{t122} 


{
\begin{proposition}\label{t3} 
	$\ $
 \begin{enumerate}[i)]
     \item With  function $\!\Phi\!$ in (\ref{potential-fun}) and payoffs $J_{i}$ for $i\!\in\!\mathcal{N}_s$ in \eqref{f1},     game ${{G}}\!=\!\{\mathcal{N}_{s},\! \{\Omega_i\}_{i\in\!\mathcal{N}_{s}}, \!\{J_i\}_{i\in\!\mathcal{N}_{s}}\}$ is a  potential  game.
%
 \item  {Under Assumption 1,  
 the  NE  $\boldsymbol{x}^{\Diamond}$ of  potential game ${{G}}$ is unique and is equal to the solution $\boldsymbol{x}^{\star}$  of \eqref{SNL}.
  }
 \end{enumerate}
\end{proposition}
Proposition \ref{t3} indicates the alignment between each non-anchor node’s individual goal and the whole network’s objective. 
 The  NE ensures not only  that each non-anchor node can adopt its optimal location strategy from the individual perspective, but also that  the sensor network as a whole can achieve a precise localization  from the global perspective, i.e., 
 the strategy profile  $\boldsymbol{x}^{\Diamond}$ satisfies $\|x_{i}^{\Diamond}-x_{j}^{\Diamond}\|^2-d_{i j}^2=0$ and $\|x_{i}^{\Diamond}-a_{l}\|^2-e_{i l}^2=0$ for any  $(i,j)\in \mathcal{E}_{ss}$ and $(i,l)\in \mathcal{E}_{as}$.

 
%
%
 
\begin{remark}
{
It should be noted that $J_{i}$ and $	\Phi$  are non-convex functions in both \cite{jia2013distributed,ke2017distributed} and our model, which brings the challenge to obtain the NE. However, different from the use of the Euclidean norm in \cite{jia2013distributed,ke2017distributed}, i.e., $\|\|x_{i}-x_{j}\|\!-d_{ij}\|$, 
we adopt the square of Euclidean norm to characterize  $J_{i}$ and $\Phi$, i.e., $\|\|x_{i}-x_{j}\|^2-d_{ij}^2\|^2$. These functions endowed with continuous fourth-order polynomials enable us to avoid the non-smoothness in \cite{jia2013distributed,ke2017distributed}, so as to deal with the intrinsic non-convexity of SNL with useful technologies,
which will be shown in Section \ref{sec:duality}.
In fact, previous efforts fail to adequately address the intrinsic 
non-convexity of SNL, either under potential games \cite{jia2013distributed,ke2017distributed}
or other modeling methods \cite{calafiore2010distributed,wan2019sensor}. Thus, they
merely yield an approximate solution or a local NE by relaxing non-convex constraints or relying on additional convex assumptions. 
}
\end{remark} 

On this basis, the problems of finding the global sensor network solution $\boldsymbol{x}^{\star}$ in \eqref{SNL}
is converted into seeking the global NE $\boldsymbol{x}^{\Diamond}$ of game $G$.
As for convex games, most of the existing research
works seek NE   via investigating first-order stationary points \cite{facchinei2010penalty,koshal2016distributed,chen2021distributed}.
 However, in such a non-convex regime (\ref{potential-fun}), 
one cannot expect to find an NE easily following this way, because a stationary point in non-convex settings is not equivalent to an  NE anymore. 
 The seeking of
 the global NE have remained ambiguous in the context of this non-convex SNL problem. 
{
 We expect to solve the following problems:
 \begin{problem}
 $~$
 \begin{enumerate}[(i)]
 	\item Under what conditions, a stationary point of (\ref{potential-fun}) is consistent to  the  NE $\boldsymbol{x}^{\Diamond}$?
 	\item How to design an efficient algorithm to seek the  NE $\boldsymbol{x}^{\Diamond}$?
 \end{enumerate}
  \end{problem}
 To this end, we solve the above questions by following processes. First, we establish the relationship between the NE and the stationary points by virtue of canonical dual theory and investigate the necessary and sufficient condition for reaching
 the  NE  (in Section \ref{sec:process}).  Then, we design a conjugate-based algorithm with the assisted complementary dual information, and establish the global convergence to the NE (in Section \ref{sec:cent}).  Finally, we consider a distributed implementation of the NE-seeking algorithm to reduce the  transmission of dual information
in the network 
(in Section \ref{sec:dis}).
%


 \section{Derivation of global Nash equilibrium}\label{sec:duality}

{
In this section, we   explore 
the  
necessary and sufficient
condition for a stationary point to coincide with the NE
and develop a conjugate-based SNL algorithm to compute the NE. }  


\subsection{Dealing with non-convexity}\label{sec:process}
 It is hard to directly identify and verify  whether a stationary point is the  NE
  on  the non-convex potential function (\ref{potential-fun}).  Here, 
we employ canonical duality theory \cite{gao2017canonical} to transform (\ref{potential-fun}) into a complementary function and  investigate the relationship between  a stationary point of the dual problem and the  NE of game (\ref{f1})
with the following two steps.

\noindent\textbf{Step 1}
We first reformulate  (\ref{potential-fun}) in canonical form. Consider the vectors
  \begin{equation}\label{ci} 
  	\begin{aligned}
  &{\xi}^s={\psi}^s(\boldsymbol{x})=\operatorname{col}\{\xi_{ij}^{s}\}_{(i,j)\in\mathcal{E}_{ss} }\in \Theta_s,\\
  &{\xi}^a={\psi}^a(\boldsymbol{x})=\operatorname{col}\{\xi_{il}^{a}\}_{(i,l)\in\mathcal{E}_{as} }\in \Theta_a,
  \end{aligned}
\end{equation}
where $\xi_{ij}^{s}=\|x_{i}-x_{j}\|^2$ and $\xi_{il}^{a}=\|x_{i}-a_{l}\|^2$. Define the quadratic functions
\begin{equation}\label{cl}
\begin{aligned}
	& \Psi_{s}({\xi}^s)\!= \!\!\!\!\!\!\!\sum_{(i, j) \in \mathcal{E}_{ss}}\!(\xi_{i j}^s\!-\!d_{i j}^2)^2,\; \Psi_{a}({\xi}^a)\!=\!\!\!\!\!\!\sum_{(i, l) \in \mathcal{E}_{as}}\!(\xi_{i l}^a\!-\!e_{i l}^2)^2.
\end{aligned}
\end{equation}
Then, the potential function (\ref{potential-fun}) can  be  rewritten as:
$$
\Phi(\boldsymbol{x})=\Psi_{s}({\psi}^s(\boldsymbol{x}))+\Psi_{a}({\psi}^a(\boldsymbol{x})). 
$$
{ Moreover, it follows from  \cite{gao2017canonical} that 
	both $ \Psi_{s}: \Theta_{s} \rightarrow \mathbb{R}$  and $ \Psi_{s}: \Theta_{s} \rightarrow \mathbb{R}$
are  convex differential
canonical functions, where the derivatives $\nabla\Psi_{s}({\xi}^s): \Theta_{s}\rightarrow \Theta^{*}_{s}$ and $\nabla\Psi_{a}({\xi}^a): \Theta_{a}\rightarrow \Theta^{*}_{a}$
	are one-to-one mappings from $\Theta_{s}$ and $\Theta_{a}$ to their ranges $\Theta^{*}_{s}$ and  $\Theta^{*}_{a}$.
 }
This indicates that the Legendre conjugates of $\Psi_{s}$ and $\Psi_{a}$ can be uniquely defined by 
 \begin{equation}\label{cr}
 	\begin{aligned}
 		\Psi_{s}^*({\sigma}^s) & = (\xi^s)^{\operatorname{T}} \sigma^s-\Psi_{s}(\xi^s)=
 		\sum\nolimits_{(i, j) \in \mathcal{E}_{ss}} \frac{1}{4}(\sigma_{i j}^s)^2+d_{i j}^2 \sigma_{i j}^s, \\
 		\Psi_{a}^*({\sigma}^a) & =(\xi^a)^{\operatorname{T}} \sigma^a-\Psi_{a}(\xi^a)=\sum\nolimits_{(i, l) \in \mathcal{E}_{as}} \frac{1}{4}(\sigma_{i l}^a)^2+e_{i l}^2 \sigma_{i l}^a,
 	\end{aligned}
 \end{equation}
where \begin{equation*}
	\begin{array}{ll}
		\sigma_{i j}^s=\nabla_{\xi_{i j}^s}\Psi_{s}({\xi}^s)=2(\xi_{i j}^s-d_{i j}^2), & (i, j) \in \mathcal{E}_{ss}, \\
		\sigma_{i l}^a=\nabla_{\xi_{i l}^a}\Psi_{a}({\xi}^a)=2(\xi_{i l}^e-e_{i l}^2), & (i, l) \in \mathcal{E}_{as}, \vspace{-0.2cm}
	\end{array}
\end{equation*}
with
$$\sigma^{s}=\operatorname{col}\{\sigma_{ij}^{s}\}_{(i,j)\in\mathcal{E}_{ss} }\in\Theta^{*}_{s} ,\quad\sigma^{a}=\operatorname{col}\{\sigma_{il}^{a}\}_{(i,l)\in\mathcal{E}_{as} }\in\Theta^{*}_{a}.$$
Also, the following  duality relations are invertible on  both $\Theta_s \times \Theta^*_s$ and $\Theta_a \times \Theta^*_a$:
\begin{equation}\label{dual}
	\begin{array}{ll}
		\sigma_{i j}^s\!=\!\nabla_{\xi_{i j}^s}\Psi_{s}({\xi}^s)\! \iff \!\xi_{i j}^s\!=\!\nabla_{\sigma_{i j}^s}\Psi^*_{s}({\sigma}^s), & \!\!\!\!(i, j)\! \in \!\mathcal{E}_{ss}, \\
		\sigma_{i l}^a\!=\!\nabla_{\xi_{i l}^a}\Psi_{a}({\xi}^a)\!\iff\! \xi_{i j}^a\!=\!\nabla_{\sigma_{i j}^a}\Psi^*_{a}({\sigma}^a), & \!\!\!\!(i, l) \!\in\! \mathcal{E}_{as}. \vspace{-0.1cm}
	\end{array}
\end{equation}
For convenience,  define the compact form $$\sigma=\operatorname{col}\{\sigma^{s},\sigma^{a} \},\quad \sigma\in\Theta^* =\Theta_s^*\times\Theta_a^*\subseteq \mathbb{R}^{q},$$ 
where {$q=\left|\mathcal{E}_{ss}\right|+\left|\mathcal{E}_{as}\right|$ is the total number of elements in the edge sets $\mathcal{E}_{ss}$ and $\mathcal{E}_{as}$.} We regard  $\sigma$ as a canonical dual variable on the whole dual space. Then
we define the following complementary function  $\Gamma: \boldsymbol{\Omega}\times \Theta^{*}\rightarrow \mathbb{R} $,
	\begin{align}\label{complementary}
		\Gamma\left(
  \boldsymbol{x},
  \sigma\right)=&(\xi^s)^{\operatorname{T}} \sigma^s-\Psi_{s}^*(\sigma^s)+(\xi^a)^{\operatorname{T}} \sigma^a-\Psi_{a}^*(\sigma^a)\nonumber\\
		=&\!\!\sum\nolimits_{(i,j)\in\mathcal{E}_{ss} }\!\!\!\!\sigma_{ij}^{s}(\|x_{i}-x_{j}\|^{2}-d_{ij}^2)\!\!\nonumber\\
		&\!\!+\!\!\sum\nolimits_{(i,l)\in\mathcal{E}_{as} }\!\!\!\!\sigma_{il}^{a}(\|x_{i}-a_{l}\|^{2}-e_{il}^2)\\
		&\!\!-\!\!\sum\nolimits_{(i,j)\in\mathcal{E}_{ss} }\!\!\!\!\frac{(\sigma_{ij}^{s})^2}{4}\!\!-\!\!\sum\nolimits_{(i,l)\in\mathcal{E}_{as} }\!\!\!\!\frac{(\sigma_{il}^{a})^2}{4}. \nonumber
	\end{align}
So far, we have transformed the non-convex function (\ref{potential-fun}) into the complementary function (\ref{complementary}). It can be established 
that
 if 	$(\boldsymbol{x}^{\Diamond},{\sigma}^{\Diamond})  \in \boldsymbol{\Omega}\times\Theta^{*} $ 
is a stationary point of    $ \Gamma (\boldsymbol{x},\sigma)$, then $\boldsymbol{x}^{\Diamond}$ is a stationary point of $\Phi\left(\boldsymbol{x}\right)$. This equivalency is because of the fact that
 the duality relations (\ref{dual}) are  unique and invertible.

\noindent\textbf{Step 2}
 We introduce a sufficient feasible domain for the introduced conjugate variable $\sigma$, in order to  investigate the global optimality of the stationary points in (\ref{complementary}).
Consider the  second-order derivative
of $\Gamma\left(\boldsymbol{x},\sigma\right)$ in $\boldsymbol{x}$. Due to the expression of (\ref{complementary}), we can find that $\Gamma$ is quadratic in $\boldsymbol{x}$. Thus, $\nabla^{2}_{\boldsymbol{x}}\Gamma$ is   $\boldsymbol{x}$-free,  and is indeed a linear combination for the elements of $\sigma$.
In this view, we denote $	P(\sigma)= \nabla^{2}_{\boldsymbol{x}}\Gamma$. The specific expression of $\nabla^{2}_{\boldsymbol{x}}\Gamma$ is given in Appendix \ref{aa}.  On this basis, we introduce the following set of $\sigma$ 
\begin{equation}\label{s23}
 \mathbb{E}^{+}= \Theta^{*}\cap \{\sigma : 
	P(\sigma) \succeq \boldsymbol{0}_{nN}
	\}. 
\end{equation} 
$\mathbb{E}^{+}$ can be regarded as a sufficient convex feasible domain for $\sigma$, and is proved to be nonempty (see in Appendix \ref{aa}).  When $ \sigma\in\mathbb{E}^{+}  $,  the positive semidefiniteness of $P(\sigma)$ implies that $\Gamma(\boldsymbol{x},{\sigma}  )$ is convex with respect to $\boldsymbol{x}$. Besides, the convexity of $\Psi(\xi)$ derives that its Legendre conjugate $\Psi^{*}(\sigma)$ is also convex, implying that the  complementary function  $\Gamma(\boldsymbol{x},{\sigma}  )$  is concave in $\sigma$. {
	This convex-concave property of $\Gamma$ enables us to further investigate the identification and verification conditions of the NE.}


Recall that the  NE
can be represented as  a strategy profile  $\boldsymbol{x}^{\Diamond}$ that satisfies $\|x_{i}^{\Diamond}-x_{j}^{\Diamond}\|^2-d_{i j}^2=0$ and $\|x_{i}^{\Diamond}-a_{l}\|^2-e_{i l}^2=0$ for any $(i,j)\in \mathcal{E}_{ss}$ and $(i,l)\in \mathcal{E}_{as}$.
It can be deduced from  (\ref{dual}) that
\begin{equation}\label{zero}
	\begin{array}{ll}
		\sigma_{i j}^{s}=2(\|x_{i}^{\Diamond}-x_{j}^{\Diamond}\|^2-d_{i j}^2)=0, & \forall (i, j) \in \mathcal{E}_{ss}, \\
		\sigma_{i l}^{a}=2(\|x_{i}^{\Diamond}-a_{l}\|^2-e_{i l}^2)=0, &\forall (i, l) \in \mathcal{E}_{as},
	\end{array}
\end{equation}
These indicate that the vector of dual variables ${\sigma}^{\Diamond}$ corresponding to
the  NE ${\boldsymbol{x}^{\Diamond}}$ is such that 
$${\sigma}^{\Diamond}=\boldsymbol{0}_{q}.$$
Moreover, we can verify  $\boldsymbol{0}_{q}\in \mathbb{E}^{+}$.
Hence, based on the above two steps,  we provide the following result to connect the global NE and the stationary point in 	$\Gamma\left(\boldsymbol{x},\sigma\right)$, whose proof is in Appendix \ref{bb}.
{
	\begin{mythm}\label{t2}
	Under Assumption 1,	
	a profile $\boldsymbol{x}^{\Diamond}$ is the global NE		of  non-convex  game ${{G}}$ if and only if there exists ${\sigma}^{\Diamond}\in \mathbb{E}^{+}$ such that
	$(\boldsymbol{x}^{\Diamond},{\sigma}^{\Diamond}) \in \boldsymbol{\Omega}\times\mathbb{E}^{+} $ 
	is a stationary point of    $ \Gamma (\boldsymbol{x},\sigma)$  satisfying 
	$$\sigma_{i j}^{s\Diamond}=\nabla_{\xi_{i j}^s}\Psi_{s}({\xi}^s)|_{\xi_{i j}^s=\|x_{i}^{\Diamond}-x_{j}^{\Diamond}\|^2}, \forall (i, j) \in \mathcal{E}_{ss}, $$  
	$$\sigma_{i l}^{a\Diamond}=\nabla_{\xi_{i l}^a}\Psi_{a}({\xi}^a)|_{\xi_{i j}^a=\|x_{i}^{\Diamond}-a_{l}\|^2}, \forall (i, l) \in \mathcal{E}_{as}.$$
\end{mythm}}

The result in Theorem \ref{t2} establishes the isomorphic mapping between original potential game ${{G}}$ and complementary dual problem (\ref{complementary}) on   $\boldsymbol{\Omega}\times\mathbb{E}^{+}$.
From the sufficiency perspective,  
 since $\mathbb{E}^{+}$ is nonempty,  the canonical duality
theory can be guaranteed to solve such non-convex problems. The identification of the  NE, indeed the global NE, can be achieved by obtaining the stationary point of $ \Gamma (\boldsymbol{x},\sigma)$ on $\boldsymbol{\Omega}\times\mathbb{E}^{+}$ and checking whether the duality relations 	$\sigma_{i j}^{s\Diamond}=\nabla_{\xi_{i j}^s}\Psi_{s}({\xi}^s)|_{\xi_{i j}^s=\|x_{i}^{\Diamond}-x_{j}^{\Diamond}\|^2}$ and 
$\sigma_{i l}^{a\Diamond}=\nabla_{\xi_{i l}^a}\Psi_{a}({\xi}^a)|_{\xi_{i l}^a=\|x_{i}^{\Diamond}-a_{l}\|^2}$ are satisfied.
 From the necessity perspective, 
 once the existence of the NE is guaranteed, the validity of the duality relations 
  can be verified on $\boldsymbol{\Omega}\times\mathbb{E}^{+}$. Therefore, the global NE can be obtained by solving  $\Gamma(\boldsymbol{x},\sigma)$ via
  its first-order conditions and employing the duality relations as a criterion.
%

Additionally, since the optimal ${\sigma}^{\Diamond}=\boldsymbol{0}_{q}\in \mathbb{E}^{+}$, where $q=\left|\mathcal{E}_{ss}\right|+\left|\mathcal{E}_{as}\right|$,  we can 
replace the $\mathbb{E}^{+}$ with simple unit square constraints in the practical implementation, whose proof is shown in  Appendix \ref{cc}.
	\begin{corollary}\label{c1}
		{
 The conclusion of  Theorem \ref{t2} still holds when
 replacing $\mathbb{E}^{+}$ with $[0,W]^{q}$, where  $W$ is a positive constant. 
  }
	\end{corollary}	

\begin{remark}
Compared with the related works on solving the original non-convex SNL problems \cite{calafiore2010distributed,jia2013distributed}, we employ an auxiliary canonical dual variable $\sigma$ and construct the complementary function $\Gamma\left( \boldsymbol{x},\sigma\right)$ to handle the non-convexity.   Furthermore, we introduce an obtainable domain set $\mathbb{E}^{+}$ of $\sigma$ to identify the global solution,  which avoids  the possibility of the algorithm being stuck into stationary points or local optima,
as discussed in \cite{calafiore2010distributed,jia2013distributed}. 
 Meanwhile, we do not require the  relaxation of the non-convex rank constraints or the addition assumptions such as the good priori initialization \cite{calafiore2010distributed}  or the convex
 network structure \cite{wang2006further,srirangarajan2008distributed}. Besides, different from \cite{wan2019sensor}, we effectively obtain an exact solution with a theoretical guarantee.
	\end{remark}

\subsection{Algorithm design and analysis}\label{sec:cent}
Based on  Theorem \ref{t2}, 
we design a centralized conjugate-based algorithm (Algorithm 1) for the SNL problem with the assisted complementary information (the   Legendre conjugate of $\Psi$  and the canonical conjugate variable $\sigma$) and 
projection operators $\Pi_{\Omega_i}(\cdot)$, $\Pi_{\mathbb{E}^{+}}(\cdot)$.  
Denote  $
\boldsymbol{\zeta}=\{\boldsymbol{x},\sigma\}
$,
$
\boldsymbol{\Xi}=\boldsymbol{\Omega} \times \mathbb{E}^{+}
$, and 
	\begin{equation*}
	\begin{aligned}
		\mathrm{F}(\boldsymbol{\zeta})&\triangleq \operatorname{col}\{ \operatorname{col}\{ \nabla_{x_{i}} \Gamma\left( \boldsymbol{x},\sigma\right) \}_{i=1}^{N}, 	\nabla_{\sigma} \Gamma\left( \boldsymbol{x},\sigma\right) \}.
	\end{aligned}
\end{equation*}
\begin{algorithm}[t]
	\caption{Conjugate-based SNL Algorithm}
	\label{centralized-alg}
	\begin{algorithmic}[1]
 \renewcommand{\algorithmicrequire}{ \textbf{Input:}}
		\REQUIRE Set stepsize  $ \{\alpha[k] \} $.
  \renewcommand{\algorithmicrequire}{ \textbf{Output:}}
 \REQUIRE  $\boldsymbol{x}[k],\sigma[k]$
  \renewcommand{\algorithmicrequire}{ 
  \textbf{Initialize:}}
		\REQUIRE Set $\sigma[0]\in {\mathbb{E}^{+}}, x_{i}[0]\in \Omega_{i}, \,i\in\{1,\dots,N\}  $, 
		\FOR{$k = 1,2,\dots$}
		\STATE  update the shared  canonical dual variable:\;\vspace{0.1cm}
		
		$\sigma[k+1] = \Pi_{\mathbb{E}^{+}}(\sigma[k] +\alpha[k]
		\nabla_{\sigma} \Gamma( \boldsymbol{x}[k],\sigma[k])	
		)
		$
		\FOR{$i = 1,\dots N$}
		\STATE update the decision  variable of non-anchor node $i$:\;\vspace{0.1cm}
		
		$	x_{i}[k+1] =\Pi_{\Omega_{i}}( x_i[k]-\alpha[k] \nabla_{x_{i}} \Gamma( \boldsymbol{x}[k],\sigma[k]) )
		$
		\ENDFOR
		\ENDFOR
	\end{algorithmic}
 
\end{algorithm}
Algorithm \ref{centralized-alg}  can be rewritten in a compact form as
\begin{equation}\label{cen-compact}
\boldsymbol{\zeta}^{k+1}=\Pi_{\boldsymbol{\Xi}}[\boldsymbol{\zeta}^{k}-\alpha_{k} \mathrm{F}({\boldsymbol{\zeta}}^{k})].
\end{equation}
Consider the weighted averaged iterates in course of $k$ iterates as
$\hat{\boldsymbol{x}}[k]= {\sum_{j=1}^{k} \alpha[j] \boldsymbol{x}[j]}/{\sum_{j=1}^{k} \alpha[j]},
	\quad \hat{{\sigma}}[k]= {\sum_{j=1}^{k} \alpha[j]\sigma[j]}/{\sum_{j=1}^{k} \alpha[j]}.$ Due to the nonexpansiveness of projection operators and the convex-concave characterization of $ \Gamma (\boldsymbol{x},\sigma)$ on $\boldsymbol{\Omega} \times \mathbb{E}^{+}$, we see that Algorithm 1 globally converges 
 in the following theorem.




\begin{mythm}\label{cent}
	$\ $
 \begin{enumerate}[i)]
     \item Under Assumption 1,	a fixed point  $(\boldsymbol{x}^{\Diamond},\sigma^{\Diamond}) $ of Algorithm 1 corresponds to the  NE of game ${{G}}$ 
 if and only if 
	\begin{equation}\label{relation}
		\begin{array}{ll}
			\sigma_{i j}^{s\Diamond}=2(\|x_{i}^{\Diamond}-x_{j}^{\Diamond}\|^2-d_{i j}^2), & \forall (i, j) \in \mathcal{E}_{ss}, \\
			\sigma_{i l}^{a\Diamond}=2(\|x_{i}^{\Diamond}-a_{l}\|^2-e_{i l}^2), &\forall (i, l) \in \mathcal{E}_{as}.
		\end{array}
	\end{equation}
 \item Under Assumption 1,	Algorithm \ref{centralized-alg} converges at a $ \mathcal{O}(1/\sqrt{k}) $  rate with  step size $\alpha_k=2{{M}_2}^{-1}\sqrt{{{d}}\,/\,{2k}} $,  i.e.,
	\begin{align*}
		\Gamma(\hat{\boldsymbol{x}}[k],\sigma^{\Diamond})-
		\Gamma(\boldsymbol{x}^{\Diamond},\hat{\sigma}[k])
		\leq 
		\frac{1}{\sqrt{k}} \sqrt{2d}{M}_2,
	\end{align*}
	where ${\mu}=\min\{{\mu_{x}}/{2},{\mu_{\sigma}}/{2} \}$,  
	and ${d}$, ${M}_2$ are two positive constants. 
 \end{enumerate}
\end{mythm}
The proof of Theorem \ref{cent} can be found in Appendix \ref{bbw}.
 The convergence of Theorem \ref{cent} is described by the duality gap within the complementary function, where we prove that Algorithm 1 can converge to the NE with a satisfied duality relation \eqref{dual}. 


	In fact, it is necessary to  employ the duality relation  as a criterion to check or verify 
 the convergent point $(\boldsymbol{x}^{\Diamond},\sigma^{\Diamond})$ of Algorithm 1, because the computation of   $\sigma^{\Diamond}$ is restricted on the sufficient domain $\mathbb{E}^{+}$ instead of the original $\Theta^*$. In this view, 	
	the gradient of $\sigma^{\Diamond}$ may fall into the normal cone $\mathcal{N}_{\mathbb{E}^{+}}(\sigma^{\Diamond})$ instead of being equal to $\boldsymbol{0}_{q}$, thereby losing the one-to-one relationship with $\boldsymbol{x}^{\Diamond}$. Thus,  $\boldsymbol{x}^{\Diamond}$ may not be the global NE.
	In addition, we cannot directly employ the standard Lagrange multiplier
	method and the associated Karush-Kuhn-Tucker (KKT) theory herein, because we need	  to first confirm a feasible  domain of $\sigma$ by utilizing
	canonical duality information (referring to $\Theta^*$). In other words, once the duality relation is verified, we can say that the  convergent point $(\boldsymbol{x}^{\Diamond},\sigma^{\Diamond})$ of Algorithm 1 is indeed the global NE of game \eqref{f1}.

We summarize a road map for seeking 
	NE in this  SNL problem for friendly comprehension. That is, 
	once	the problem is defined and formulated,  we first transform the original non-convex potential function into a  complementary dual function. Then we seek the stationary point of  $ \Gamma (\boldsymbol{x},\sigma)$  via algorithm iterations, wherein the dual variable $\sigma$
	is restricted on  $\mathbb{E}^{+}$. Finally, after
	obtaining the stationary point by convergence, we
 identify whether the convergent point satisfies the
	duality relation. If so, the
	convergent point is the  NE (global NE).

\begin{remark}
	Canonical duality theory  has also been  employed in some classic optimization works to solve non-convex problems \cite{zhu2012approximate,latorre2016canonical,ren2021distributed}.  
	Note that the non-convex optimization problems considered in these works are different from the SNL game problem in this paper. Also, unlike  the relaxation approach adopted in \cite{zhu2012approximate}  to obtain approximate solutions, our goal is to compute an exact global solution to the original problem.  Moreover, the set $\mathbb{E}^{+}$ is guaranteed to be non-empty in this SNL problem,  so that we can make good use of duality theory and	
	establish an isomorphic mapping between the global solution and the stationary point of the complementary dual problem.  This is different from the sufficient conditions provided in \cite{latorre2016canonical}  and  \cite{ren2021distributed}, which are based on the imposed assumption of nonempty $\mathbb{E}^{+}$.
	
	
\end{remark} 


\section{Distributed implementation}\label{sec:dis}


{In this section, we consider 
a distributed implementation for the NE-seeking algorithm.

Note that  set constraints $\Omega_i$ and decision variables $x_i$ are non-anchor nodes' private information. Moreover,
the implementation of Algorithm 1 in reality relates to the amount of information exchange required in the networks, because each non-anchor node needs to know the global information of the canonical dual variable $\sigma$ when calculating the gradient $\nabla_{x_{i}} \Gamma( \boldsymbol{x},\sigma) $. This is often hard to accomplish in reality. To improve the scalability for solving large-scale SNL problems, 
 our goal here is to reduce the information exchange requirements by taking Algorithm 1 into a distributed implementation. 

To this end,  we propose a symmetrically distributed implementation with the help of sliding mode control \cite{sarpturk1987stability} and extra-gradient methods \cite{korpelevich1976extragradient} in Section \ref{disa}.  Then we provide explanations about the
algorithm mechanism and design in Section  \ref{dis-b}, and finally give the convergence analysis of the distributed algorithm in Section  \ref{dis-c}.}
\vspace{-0.3cm}
\subsection{Distributed algorithm design}\label{disa}
Consider a distributed situation that 
each player only accesses the information of its neighbors, and exchanges information with them. 
For $i\in \mathcal{N}_s$,
if non-anchor node $j$ is a neighbor of non-anchor node $i$, 
or in other words
$(i,j)\in\mathcal{E}_{ss}$, then the  element $\sigma_{ij}^{s} \in \mathbb{R}$ in $\sigma$   is only associated with non-anchor nodes $i$ and  $j$. Thus, non-anchor nodes $i$ and $j$ need to control their local variables  $\sigma_{ij}^{s_{i}}\in \mathbb{R}$ and $\sigma_{ij}^{s_{j}}\in \mathbb{R}$ to estimate $\sigma_{ij}^{s}$,  respectively,   while others not. On the other hand, if anchor $k$ is a neighbor of non-anchor node $i$, or
 in other words $(i,l)\in\mathcal{E}_{as}$, $\sigma_{il}^{a}\in \mathbb{R}$ in $\sigma$  is only associated with non-anchor node $i$ itself since the location of anchor $k$ is known. Thus, non-anchor node $i$  controls a local variable $\sigma_{il}^{a_{i}}\in \mathbb{R}$ to calculate $\sigma_{il}^{a}$ individually.   


  Denote the local canonical dual variable profile for non-anchor node $i$ as
\begin{equation}\label{sigma1}
 \sigma_{i}=\operatorname{col}\{\operatorname{col}\{\sigma_{ij}^{s_{i}}\}_{j\in\mathcal{N}_{s}^{i} },\operatorname{col}\{\sigma_{il}^{a_{i}}\}_{l\in\mathcal{N}_{a}^{i} }\}\subseteq \mathbb{R}^{q_i},
 \end{equation}
 where $q_i=|\mathcal{N}_{s}^{i}|+|\mathcal{N}_{a}^{i}| $ is the total number of elements in the $i$'s neighbor sets $\mathcal{N}_{s}^{i}$ and $\mathcal{N}_{a}^{i}$.  
Recalling  Corollary \ref{c1}, we can take the  constrains of both $\sigma_{ij}^{s_{i}}$ and $\sigma_{il}^{a_{i}}$ as  $[0,W]$. 
In this way, we denote the feasible set of $\sigma_{i}$ by
$$\mathscr{E}_{ i}^{+}=\prod\nolimits_{j\in\mathcal{N}_{s}^{i}}\mathscr{E}_{ ij}^{s+}\times\prod\nolimits_{l\in\mathcal{N}_{a}^{i}}\mathscr{E}_{ il}^{a+}= [0,W]^{q_i}, $$   where $\sigma_{ij}^{s_{i}}\in \mathscr{E}_{ ij}^{s+}=[0,W]$ and  $\sigma_{il}^{a_{i}}\in \mathscr{E}_{ il}^{a+}=[0,W]$.

Thus, the global complementary function $\Gamma(\boldsymbol{x},\sigma)$ in \eqref{complementary} is decomposed into a group of local functions, i.e., for $i\in\mathcal{N}_s$, 
\begin{equation}\label{dis-comle}
\begin{aligned}
&\varGamma_{i}\!(x_{i},\{x_j\}_{j\in\mathcal{N}_{s}^{i} },\!\sigma_i,\!\{\sigma_j\}_{j\!\in\!\mathcal{N}_{s}^{i} })\!=\!\!\!\sum_{j\in\mathcal{N}_{s}^{i}  }\!\!\!(\!\frac{\sigma_{ij}^{s_{i}}\!+\!\sigma_{ij}^{s_{j}}}{2})\!(\|x_{i}\!-\!x_{j}\|^{2}\!\!-\!d_{ij}^2\!)\\
	&+\!\sum_{l\in\mathcal{N}_{a}^{i}  }\sigma_{ij}^{a_{i}}(\|x_{i}\!-\!a_{l}\|^{2}-e_{il}^2)\!-\!\sum_{j\in\mathcal{N}_{s}^{i}  }\frac{(\frac{\sigma_{ij}^{s_{i}}\!+\!\sigma_{ij}^{s_{j}}}{2})^2}{4}\!-\!\sum_{l\in\mathcal{N}_{a}^{i}}\frac{(\sigma_{il}^{a_{i}})^2}{4}. \vspace{-0.2cm}
\end{aligned}
\end{equation}
Then  we propose a distributed algorithm based on  symmetric design and extra-gradient (DSDEG) in Algorithm 2
to 
solve the SNL problem, followed by design explanations in \ref{dis-b}.

Denote $ \boldsymbol{q}=\sum_{i=1}^{N}q_i$, and we obtain
$$
\bm{\mathscr{E}^{+}}=\prod_{i=1}^{N}\mathscr{E}_{ i}^{+} \subseteq \mathbb{R}^{\boldsymbol{q}}, \quad \boldsymbol{\sigma}= \operatorname{col}\{\sigma_{1}, \dots, \sigma_{N}\}  \in{\mathscr{E}}^{+} \subseteq \mathbb{R}^{\boldsymbol{q}}.\vspace{-0.3cm}
$$
\begin{algorithm}[t]
	\caption{DSDEG-based SNL Algorithm}
	\label{min-max}
	\begin{algorithmic}[1]
\renewcommand{\algorithmicrequire}{ \textbf{Input:}}
		\REQUIRE Set stepsize  $ \beta>0$.
  \renewcommand{\algorithmicrequire}{ \textbf{Output:}}
  \REQUIRE 
 $\left\{x_i[k]\right\}_{k=1}^{\infty}, \left\{\sigma_i[k]\right\}_{k=1}^{\infty}$.
\renewcommand{\algorithmicrequire}{ \textbf{Initialize:}}
  \REQUIRE Set {$x_{i} \![0]$, $ \!\tilde{x}_{i}[0]\!\!\in\! \mathbb{R}^{n_{i}}$, $ \!{\sigma}_{ij}^{s_{i}} \![0]$,$ \tilde{\sigma}_{ij}^{s_{i}} \![0],\!{\sigma}_{il}^{a_{i}}[0],\! \tilde{\sigma}_{il}^{a_{i}}[0]\!\in\! \mathbb{R}$ with $ {\sigma_{ij}^{s_{i}}}[0]=\sigma_{ij}^{s_{j}}[0]$ for $(i,j)\in \mathcal{E}_{ss}$.}
		\FOR{$k = 1,2,\dots$}
		\FOR{$i\in\mathcal{N}_{s}$}
			\STATE {set waiting state}: \\
			 $\tilde{x}_i[k]:= x_i[k]$, $\tilde{\sigma}_{ij}^{s_{i}}[k]:={\sigma}_{ij}^{s_{i}}[k]$, $\tilde{\sigma}_{il}^{a_{i}}[k]:={\sigma}_{il}^{a_{i}}[k]$\\
\STATE \textbf{for} $\ell = 1,2$ \textbf{do}\\
	\STATE 	\quad update decision  variable:\\
		\quad$	\tilde{x}_i[k]\! :=\!\Pi_{\Omega_{i}}( x_i[k]\!-\!\beta (\! \sum\limits_{j\!\in\!\mathcal{N}_{s}^{i} }\!(\tilde{\sigma}_{ij}^{s_{i}}[k]\!+\tilde{\sigma}_{ij}^{s_{j}}[k])(\tilde{x}_{i}[k]-$\\
		\quad\quad\quad\quad\quad\quad\;$\tilde{x}_{j}[k])\!+\!\!\sum\limits_{l\in\mathcal{N}_{a}^{i} }\!2\tilde{\sigma}_{il}^{a_{i}}\!(\tilde{x}_{i}[k]\!-\!a_{l})))$\\		
		\STATE\quad update dual variable  (non-anchor neighbor) \\ 
             \STATE\quad \textbf{for} {$ j\in\mathcal{N}_{s}^{i}$}  \textbf{do}\\
		\quad$\tilde{\sigma}_{ij}^{s_{i}}[k]\!:=\!\Pi_{\mathscr{E}_{ ij}^{s+}}( \!\sigma_{ij}^{s_{i}}[k]\!+\!\beta (\frac{1}{2}(\|\tilde{x}_{i}[k]\!-\!\tilde{x}_{j}[k]\|^{2}\!-\!d_{ij}^2)-$\\
		\quad\quad\quad\quad\quad\quad\quad\,$\!\frac{1}{8}(\tilde{\sigma}_{ij}^{s_{i}}[k]+\tilde{\sigma}_{ij}^{s_{j}}[k])),$\\
   \STATE\quad \textbf{end for}
 \STATE\quad update dual variable  (anchor neighbor) \\ 
    \STATE\quad \textbf{for} {$ l\in\mathcal{N}_{a}^{i}$}  \textbf{do}\\
	 	\quad$\tilde{\sigma}_{il}^{a_{i}}[k]\!:=\!\Pi_{\mathscr{E}_{ il}^{a+}}( \sigma_{il}^{a_{i}}\![k]\!+\!\beta(\|\tilde{x}_{i}[k]\!-\!a_{l}\|^{2}\!-\!e_{il}^2\!-\!\frac{1}{2}\tilde{\sigma}_{il}^{a_{i}}\![k])$\\ 
 \STATE\quad \textbf{end for}
 \STATE \textbf{end for}
	\STATE {set new state}: \\
$ x_i[k\!+\!1]\!\!:=\!\tilde{x}_i[k]$, ${\sigma}_{ij}^{s_{i}}[k\!+\!1]\!\!:=\!\tilde{\sigma}_{ij}^{s_{i}}[k]$, ${\sigma}_{il}^{a_{i}}[k\!+\!1]\!\!:=\!\tilde{\sigma}_{il}^{a_{i}}[k]$\\
		\ENDFOR
		\ENDFOR
	\end{algorithmic}
\end{algorithm}
 Define $
{\boldsymbol{z}}=\operatorname{col}\left\{ \boldsymbol{x},\boldsymbol{\sigma}\right\}
$, $
{\tilde{\boldsymbol{z}}}=\operatorname{col}\left\{ \tilde{\boldsymbol{x}},\tilde{\boldsymbol{\sigma}}\right\}
$ with $\tilde{\boldsymbol{x}}=\operatorname{col}\{\tilde{x_i}\}_{i=1}^N$ and $\tilde{\boldsymbol{\sigma}}=\operatorname{col}\{\operatorname{col}\{\operatorname{col}\{\tilde{\sigma}_{ij}^{s_{i}}\}_{j\in\mathcal{N}_{s}^{i} },$ $\operatorname{col}\{\tilde{\sigma}_{il}^{a_{i}}\}_{l\in\mathcal{N}_{a}^{i} }\}\}_{i=1}^N$, and $
{\boldsymbol{\varXi}}=\boldsymbol{\Omega} \times \mathscr{E}^{+}
$. 
Then the  pseudo-gradient of  \eqref{dis-comle} for $i\in \mathcal{N}_s$ can be rewritten as
\begin{align}\label{dis-F}
	\boldsymbol{F}({\boldsymbol{z}})&\triangleq
	\left[\begin{array}{c}
		\operatorname{col}\{ \nabla_{x_{i}} \varGamma_{i}(x_{i},\!\{\!x_j\!\}_{j\in\mathcal{N}_{s}^{i} },\sigma_i,\!\{\!\sigma_j\!\}_{j\in\mathcal{N}_{s}^{i} })\}_{i=1}^{N}
		\\
		\operatorname{col}\{-\nabla_{\sigma_i} \varGamma_{i}(x_{i},\!\{\!x_j\!\}_{j\in\mathcal{N}_{s}^{i} },\sigma_i,\!\{\!\sigma_j\!\}_{j\in\mathcal{N}_{s}^{i} }) \}_{i=1}^{N}
	\end{array}\right].
\end{align}
where $\nabla_{x_{i}} \varGamma_{i}=\sum_{j\in\mathcal{N}_{s}^{i} }\!(\sigma_{ij}^{s_{i}}\!\!+\sigma_{ij}^{s_{j}})(x_{i}-x_{j})\!\!+\!\!\sum_{l\in\mathcal{N}_{a}^{i} }\!\!2\sigma_{il}^{a_{i}}(x_{i}\!-\!a_{l}),$
\begin{align*}
	&\nabla_{\sigma_{i}}\varGamma_{i}=
	\left[\begin{array}{c}
		\operatorname{col}\{\frac{1}{2}(\|x_{i}-x_{j}\|^{2}\!-\!d_{ij}^2)\!-\!\frac{1}{8}(\sigma_{ij}^{s_{i}}+\sigma_{ij}^{s_{j}}) \}_{j\in\mathcal{N}_{s}^{i} }\\ \operatorname{col}\{\|x_{i}\!-\!a_{l}\|^{2}-e_{il}^2\!-\!\frac{1}{2}\sigma_{il}^{a_{i}}\}_{l\in\mathcal{N}_{a}^{i} }
	\end{array}\right].
\end{align*}
Define $\boldsymbol{F}(\tilde{\boldsymbol{z}})$ in a similar way. Algorithm 2 be rewritten in a compact form as
\begin{equation}\label{dis-compact}
	\left\{\begin{array}{l}
		\tilde{\boldsymbol{z}}[k]=\Pi_{{\boldsymbol{\varXi}}}[\boldsymbol{z}[k]-\beta \boldsymbol{F}(\boldsymbol{z}[k])], \\
		\boldsymbol{z}[k+1]=\Pi_{{\boldsymbol{\varXi}}}[\boldsymbol{z}[k]-\beta \boldsymbol{F}(\tilde{\boldsymbol{z}}[k])].
	\end{array}\right.
\end{equation} 

\subsection{Interpretation of DSDEG}\label{dis-b}


In DSDEG,  
the global canonical variable $\sigma$ is computed in a distributed way with the local information of each non-anchor node and its neighbors. 
For each neighbor $j\in\mathcal{N}_s^i$, non-anchor node $i$ 
exchange its local estimate  $\sigma_{ij}^{s_{i}}$  with $j$  and  employ the average estimation term ${(\sigma_{ij}^{s_{i}}+\sigma_{ij}^{s_{j}})}/{2}$  in  \eqref{dis-comle} to track the real $\sigma_{ij}^{s}$ in   \eqref{complementary}.  When $\sigma_{ij}^{s_{j}}=\sigma_{ij}^{s_{i}}$, $({\sigma_{ij}^{s_{i}}+\sigma_{ij}^{s_{j}}})/{2}=\sigma_{ij}^{s}$. For each neighbor $l\in\mathcal{N}_{a}^i$, the value of $\sigma_{il}^{a}$ depends on non-anchor node $i$ itself and
non-anchor node $i$ employs the local variable $\sigma_{il}^{a_{i}}$ to calculate   $\sigma_{il}^{a}$ individually.
In this regard,  non-anchor node $i$ determines its location  by taking into account its own local  information $x_i$, $\sigma_{ij}^{s_{i}}$ and $\sigma_{il}^{a_{i}}$, as well as the available observations of its neighbor $j$'s $x_j$ and $\sigma_{ij}^{s_{j}}$ for $j\in\mathcal{N}_s^i$. 
Non-anchor node $i$  calculates the local decision
variable $x_i$ and the local conjugate variable $\sigma_i$ following the  gradient descent $\nabla_{x_{i}}\varGamma_{i}$ and ascent $\nabla_{\sigma_{i}}\varGamma_{i}$ via projections. 

 
{
\begin{remark}
A major concern of the distributed implementation is the communication cost, 
which is closely related to the information assignment between each non-anchor node and its neighbors.   
The existing mainstream distributed ideas \cite{costa2006distributed,ke2017distributed,calafiore2010distributed} usually require non-anchor node $i$ to control a local estimate $\bar{\sigma_i}\in \mathbb{R}^{q}$ of the global variable $\sigma\in\mathbb{R}^{q}$ and exchange this estimate with its neighbors $\mathcal{N}_{s}^{i} $. This indicates that the communication complexity of non-anchor node $i$ at each iteration of these methods \cite{costa2006distributed,ke2017distributed,calafiore2010distributed}
is $\mathcal{O}(W^q)$.
However, note that the global variable $\sigma=\operatorname{col}\{\operatorname{col}\{\sigma_{ij}^{s}\}_{(i,j)\in\mathcal{E}_{ss}},\operatorname{col}\{\sigma_{il}^{a}\}_{(i,l)\in\mathcal{E}_{as}} \}\in\mathbb{R}^{q}$ contains the relationship among all  sensor nodes.
Thus, the above  designs may result in the local  $\bar{\sigma_i}$ containing some redundant information irrelevant to non-anchor node $i$ itself,  e.g., the information $\sigma^s_{kj}$  related to non-anchor node $j$ and non-anchor node $k$ for $(k,j)\in\mathcal{E}_{ss}$.   The advantage of using our distributed implementation over these previous works is to avoid redundant space. In fact, it follows from \eqref{sigma1} that non-anchor node $i$ does not need to save the irrelevant estimates with itself, and the communication complexity of non-anchor node $i$ at each iteration of our method is $\mathcal{O}(W^{q_i})$, instead of $\mathcal{O}(W^q)$.
\end{remark}}
 
 
 Then  we give explanations about   DSDEG's designs. 
 
{
We introduce the distributed symmetric design about $({\sigma_{ij}^{s_{i}}+\sigma_{ij}^{s_{j}}})/{2}$.
The design behind it is to guarantee the monotonicity of the pseudo-gradient $	\boldsymbol{F}({\boldsymbol{z}})$ or in other words, the positive semidefiniteness of the Jacobian matrix of  $\boldsymbol{F}({\boldsymbol{z}[k]}) $, which is important for algorithm convergence in game problems.  Denote the the Jacobian matrix of  $\boldsymbol{F}({\boldsymbol{z}[k]}) $ as $\mathcal{J}_{\boldsymbol{F}}({\boldsymbol{z}[k]}) $. 
  The design is inspired by sliding mode control \cite{sarpturk1987stability}. We give an intuitive explanation for the design idea by removing the projection $\Pi_{{\boldsymbol{\varXi}}} $ in algorithm \eqref{dis-compact} and simply considering once update, i.e., $\boldsymbol{z}[k+1]=\boldsymbol{z}[k]-\beta \boldsymbol{F}({\boldsymbol{z}}[k])$, which can be rewritten as
  $\boldsymbol{z}[k+1]=\boldsymbol{z}[k]-\beta( {\boldsymbol{F}^{\prime}}(\boldsymbol{z}[k])+\boldsymbol{u}[k])$. 
  Here, 
  ${\boldsymbol{F}^{\prime}}(\boldsymbol{z}[k])$ is the pesudo-gradient for the case where we only employ $\sigma_{ij}^{s_{i}}$ to decompose $\Gamma(\boldsymbol{x},\sigma)$, while $\boldsymbol{u}[k]$ is a state feedback controller. By  taking $\mathcal{J}_{\boldsymbol{F}}({\boldsymbol{z}[k]}) $    as a sliding manifold,
once we control $\boldsymbol{u}[k]$ to make $(\mathcal{J}_{{\boldsymbol{F}}}({\boldsymbol{z}[k]})^{\operatorname{T}}+\mathcal{J}_{{\boldsymbol{F}}}({\boldsymbol{z}[k]}))/2 $ lie with a space of  positive semidefinite matrices,
  in an ideal situation, its state variables $\boldsymbol{z}[k]$ will  ``slide'' to the convergent point without additional control. Therefore,   
  we substitute $\sigma_{ij}^{s_i}$ with the average  $({\sigma_{ij}^{s_{i}}+\sigma_{ij}^{s_{j}}})/{2}$ to decompose $\Gamma(\boldsymbol{x},\sigma)$ and 
  construct the controller as
   \begin{align}\label{contro}
 	&\boldsymbol{u}[k]={\boldsymbol{F}^{\prime}}(\boldsymbol{z}[k])-\boldsymbol{F}({\boldsymbol{z}[k]})\\
 	&\!=\!	
 	\left[\begin{array}{c}
 		\!\operatorname{col}\{\!\sum_{j\in\mathcal{N}_{s}^{i} }\!(\sigma_{ij}^{s_{i}}\!\!-\sigma_{ij}^{s_{j}})(x_{i}-x_{j}) \}_{i=1}^N\!\\ 		\operatorname{col}\!\{\!\operatorname{col}\!\{-\frac{\|x_{i}\!-\!x_{j}\|^{2}\!-\!d_{ij}^2}{2}\!+\!\frac{3\sigma_{ij}^{s_{i}}}{8}\!-\!\!\frac{\sigma_{ij}^{s_{j}}}{8}
 		\}_{j\!\in\!\mathcal{N}_{s}^{i} }, \operatorname{col}\{0\}_{l\!\in\!\mathcal{N}_{a}^{i} }\!\}_{i=1}^N \!\!\nonumber
 	\end{array}\right].
 \end{align}
 This  yields a positive
semi-definite $\mathcal{J}_{\boldsymbol{F}}({\boldsymbol{z}[k]}) $. Also, the projection on $\sigma_i$ yields that $\sigma_{ij}^{s_{i}}$ and $\sigma_{il}^{a_{i}}$ are not less than zero, which guarantees that the main diagonal elements of $\frac{1}{2}(\mathcal{J}_{{\boldsymbol{F}}}({\boldsymbol{z}[k]})^{\operatorname{T}}+\mathcal{J}_{{\boldsymbol{F}}}({\boldsymbol{z}[k]})) $  are always greater than or equal to zero. Thus, once the manifold $\mathcal{J}_{\boldsymbol{F}}({\boldsymbol{z}[k]}) $ is reached,  $\sigma_i$  are always on manifold $\mathcal{J}_{\boldsymbol{F}}({\boldsymbol{z}[k]}) $ during the iteration of Algorithm 2, which avoids chattering and high-frequency switching. 
We show more details of this design in Appendix \ref{ccf}}.

 Another design idea is the twice updates of $\boldsymbol{z}[k]$, which is to ensure the global convergence even when pseudo-gradient $\boldsymbol{F}(\boldsymbol{z})$  is merely monotone and not strictly or strongly monotone. To this end, 
the extra-gradient method 
helps us to carry
forward \cite{korpelevich1976extragradient}.   Note that  the pseudo-gradient $\boldsymbol{F}(\boldsymbol{z})$ in (\ref{dis-F}) cannot be regarded as a convex function's gradient due to the presence of  dual variables $\sigma_i$. Thus, we can not employ the  basic gradient method to make  the algorithm decrease with $\left\|\boldsymbol{z}[k]-\boldsymbol{z}^{\Diamond}\right\|^2$, which may result in an oscillation. Also, 
 the distributed SNL problem
 actually reflects a cluster of $\varGamma_i$, rather than a unified  $\Gamma$ in \eqref{complementary}. This indicates that we cannot establish the equivalency between the monotonicity of pseudo-gradient $\boldsymbol{F}(\boldsymbol{z})$ and a unified convex-concave objective function to describe the convergence.
 To overcome these difficulties, we consider slightly changing the update direction of $\boldsymbol{z}[k+1]$ 
by virtue of the extra-gradient method. 
The idea behind this technique is to compute the gradient at an (extrapolated) point different from the current point where the update is performed, that is, $$
\boldsymbol{z}[k+1]=\boldsymbol{z}[k]-\beta \boldsymbol{F}(\tilde{\boldsymbol{z}}[k]). \vspace{-0.2cm}
$$
 Also,  $\tilde{\boldsymbol{z}}[k]$ is chosen in a similar form of \eqref{dis-compact}  for the descent property.  Therefore,  the whole dynamics is guaranteed to be globally stable.

\begin{remark}
It follows from the extra-gradient methods \cite{khobotov1987modification,solodov1999new} that	the stepsize $\beta> 0$ satisfies
	\begin{equation}\label{EG-stepsize}
		\beta\cdot L <1,\vspace{-0.2cm}
	\end{equation}
	where $L>0$ is the Lipschitz constant of $\boldsymbol{F}$ over ${\boldsymbol{\varXi}}$.
	Recalling the compactness of $\boldsymbol{\Omega}$ and $\mathscr{E}^{+}$, the  Lipschitz continuity of $\nabla_{x_{i}}\varGamma_{i}$ and $\nabla_{\sigma_{i}}\varGamma_{i}$
can guarantee the existence of $L$. Intuitively, as long as $\beta $  is small enough, \eqref{EG-stepsize} can be satisfied.
In case	when obtaining  the Lipschitz constant of 	  $\boldsymbol{F}$ is challenging,
	an additional linear search procedure can be employed 
	to find a suitable stepsize, see
	\cite{solodov1999new,alakoya2021modified}.
\end{remark}

\vspace{-0.4cm}
\subsection{Algorithm analysis}\label{dis-c}

In the following, we first verify that (\ref{dis-compact}) has a fixed point that solves the distributed SNL problem.

{\begin{mythm}\label{t4}
	Under Assumption 1,	 a fixed point $(\tilde{\boldsymbol{z}}, {\boldsymbol{z}})$ of (\ref{dis-compact})  with $\tilde{\boldsymbol{z}}={\boldsymbol{z}}=\operatorname{col}(\boldsymbol{x}^{\Diamond},\boldsymbol{\sigma}^{\Diamond})$ corresponds to the NE  of game ${{G}}$ if and only if 
  \begin{equation}\label{dual-rela}
				\sigma_{i j}^{s_{i}\!\Diamond}\!\!\!=\!2(\!\|x_{i}^{\Diamond}\!\!-\!x_{j}^{\Diamond}\|^2\!\!\!-\!d_{i j}^2\!)\!\!=\!0,  \sigma_{i l}^{a_{i}\!\Diamond}\!\!\!=\!2(\!\|x_{i}^{\Diamond}\!\!-\!a_{l}\|^2\!\!\!-\!e_{i l}^2\!)\!\!=\!0. \vspace{-0.15cm}
			\end{equation}
 
\end{mythm}}

The proof of Theorem \ref{t4} can be found in Appendix \ref{c2f}.
{ Next, let $ \Upsilon $ be the set of fixed points of algorithm (\ref{dis-compact}). Since ${\boldsymbol{\varXi}}$ is compact and convex, $\boldsymbol{F}({\boldsymbol{z}})$ is continues, it follows from \cite[Theorem 1.5.5, Corollary 2.2.5]{facchinei2003finite} that
$ \Upsilon $ is nonempty. On this basis, we show the main convergence result of algorithm (\ref{dis-compact}), i.e., algorithm (\ref{dis-compact}) converges to a fixed point of $ \Upsilon $. This 
implies that the  NE can be found in a distributed way along the trajectory of (\ref{dis-compact}).}

{\begin{mythm}\label{t5}
Under Assumption 1,	
there exists ${\boldsymbol{z}}^{\Diamond}\! \!\in\! \! \Upsilon$ such that
\begin{equation}\label{lim}
	\lim _{k \rightarrow \infty} \tilde{\boldsymbol{z}}[k]=\lim _{k \rightarrow \infty} \boldsymbol{z}[k]=\boldsymbol{z}^{\Diamond}.
\end{equation}
Moreover, if 
	the convergent point $\boldsymbol{z}^{\Diamond}=(\boldsymbol{x}^{\Diamond},\boldsymbol{\sigma}^{\Diamond})$  satisfies \eqref{dual-rela},
	then $\boldsymbol{x}^{\Diamond}$ is the NE of  game ${{G}}$.
\end{mythm}}
\textbf{Proof.} 
We first prove that trajectory $\boldsymbol{z}[\cdot]$ is bounded.
Define \vspace{-0.2cm}
 \begin{equation}\label{c13}
		\boldsymbol{\upsilon}_1[k]  =\boldsymbol{z}[k]-\beta \boldsymbol{F}(\tilde{\boldsymbol{z}}[k]),\;\boldsymbol{\upsilon}_2[k]  =\boldsymbol{z}[k]-\beta \boldsymbol{F}(\boldsymbol{z}[k]).\vspace{-0.2cm}
\end{equation}
Then ${\boldsymbol{z}}[k+1]=\Pi_{\boldsymbol{\varXi}}[\boldsymbol{\upsilon}_1[k]]$ and $\tilde{\boldsymbol{z}}[k]=\Pi_{\boldsymbol{\varXi}}[\boldsymbol{\upsilon}_2[k]]$.

For any $\boldsymbol{z}^{\dag}\in \Upsilon$, define
 \begin{equation*}
	\begin{aligned}
		\Delta[k]=&\left\|\boldsymbol{z}[k+1]-\boldsymbol{z}^{\dag}\right\|^2-\left\|\boldsymbol{z}[k]-\boldsymbol{z}^{\dag}\right\|^2.\vspace{-0.2cm}
  	\end{aligned}
\end{equation*}


By (\ref{dis-compact}) and (\ref{c13}),
	\begin{align}\label{c2}
	\Delta[k]= & -\|\boldsymbol{z}[k]-\boldsymbol{z}[k+1]\|^2+2 \beta \boldsymbol{F}(\tilde{\boldsymbol{z}}[k])^\mathrm{T}\left(\boldsymbol{z}^{\dag}-\boldsymbol{z}[k+1]\right) \notag\\
		& +2 \left(\Pi_{\boldsymbol{\varXi}}[\boldsymbol{\upsilon}_1[k]]-\boldsymbol{z}^{\dag}\right)^\mathrm{T}\left(\Pi_{\boldsymbol{\varXi}}[\boldsymbol{\upsilon}_1[k]]-\boldsymbol{\upsilon}_1[k]\right)\notag\\
  =& -\|\boldsymbol{z}[k]-\tilde{\boldsymbol{z}}[k]\|^2-\|\boldsymbol{z}[k+1]-\tilde{\boldsymbol{z}}[k]\|^2 \\
  &+2 (\tilde{\boldsymbol{z}}[k]-\boldsymbol{z}[k])^\mathrm{T}(\tilde{\boldsymbol{z}}[k]-\boldsymbol{z}[k+1])\notag\\
    &+2 \beta \boldsymbol{F}(\tilde{\boldsymbol{z}}[k])^\mathrm{T}\left(\boldsymbol{z}^{\dag}-\boldsymbol{z}[k+1]\right)\notag\\
  & +2 \left(\Pi_{\boldsymbol{\varXi}}[\boldsymbol{\upsilon}_1[k]]-\boldsymbol{z}^{\dag}\right)^\mathrm{T}\left(\Pi_{\boldsymbol{\varXi}}[\boldsymbol{\upsilon}_1[k]]-\boldsymbol{\upsilon}_1[k]\right)\notag \\
 =& -\|\boldsymbol{z}[k]-\tilde{\boldsymbol{z}}[k]\|^2-\|\boldsymbol{z}[k+1]-\tilde{\boldsymbol{z}}[k]\|^2 \notag\\
&+2 \beta(\boldsymbol{z}[k+1]-\tilde{\boldsymbol{z}}[k])^\mathrm{T}(\boldsymbol{F}(\boldsymbol{z}[k])-\boldsymbol{F}(\tilde{\boldsymbol{z}}[k])) \notag\\
&+2\left(\Pi_{\boldsymbol{\varXi}}[\boldsymbol{\upsilon}_2[k]]-\boldsymbol{z}[k+1]\right)^\mathrm{T}\left(\Pi_{\boldsymbol{\varXi}}[\boldsymbol{\upsilon}_2[k]]-\boldsymbol{\upsilon}_2[k]\right)\notag\\
  &+2 \beta \boldsymbol{F}\left(\boldsymbol{z}^{\dag}\right)^\mathrm{T}\left(\boldsymbol{z}^{\dag}-\tilde{\boldsymbol{z}}[k]\right)\notag\\
  &+2 \beta \left(\boldsymbol{F}(\tilde{\boldsymbol{z}}[k])-\boldsymbol{F}\left(\boldsymbol{z}^{\dag}\right)\right)^\mathrm{T}\left(\boldsymbol{z}^{\dag}-\tilde{\boldsymbol{z}}[k]\right) \notag\\
  	& +2 \left(\Pi_{\boldsymbol{\varXi}}[\boldsymbol{\upsilon}_1[k]]-\boldsymbol{z}^{\dag}\right)^\mathrm{T}\left(\Pi_{\boldsymbol{\varXi}}[\boldsymbol{\upsilon}_1[k]]-\boldsymbol{\upsilon}_1[k]\right).\notag
	\end{align}
According to \eqref{c2}, we define 
	\begin{align*}
		H_1[k]= & -\|\boldsymbol{z}[k]-\tilde{\boldsymbol{z}}[k]\|^2-\|\boldsymbol{z}[k+1]-\tilde{\boldsymbol{z}}[k]\|^2 \notag\\
		& +2 \beta(\boldsymbol{z}[k+1]-\tilde{\boldsymbol{z}}[k])^\mathrm{T}(\boldsymbol{F}(\boldsymbol{z}[k])-\boldsymbol{F}(\tilde{\boldsymbol{z}}[k])),\\
		H_2[k]= & \left(\Pi_{\boldsymbol{\varXi}}[\boldsymbol{\upsilon}_2[k]]-\boldsymbol{z}[k+1]\right)^\mathrm{T}\left(\Pi_{\boldsymbol{\varXi}}[\boldsymbol{\upsilon}_2[k]]-\boldsymbol{\upsilon}_2[k]\right), \\
		H_3[k]= & \boldsymbol{F}\left(\boldsymbol{z}^{\dag}\right)^\mathrm{T}\left(\boldsymbol{z}^{\dag}-\tilde{\boldsymbol{z}}[k]\right), \\
		H_4[k]= & \left(\boldsymbol{F}(\tilde{\boldsymbol{z}}[k])-\boldsymbol{F}\left(\boldsymbol{z}^{\dag}\right)\right)^\mathrm{T}\left(\boldsymbol{z}^{\dag}-\tilde{\boldsymbol{z}}[k]\right), \\
		H_5[k]= & \left(\Pi_{\boldsymbol{\varXi}}[\boldsymbol{\upsilon}_1[k]]-\boldsymbol{z}^{\dag}\right)^\mathrm{T}\left(\Pi_{\boldsymbol{\varXi}}[\boldsymbol{\upsilon}_1[k]]-\boldsymbol{\upsilon}_1[k]\right),
	\end{align*}
which yields 
\begin{equation}\label{delta}
	\Delta[k]=2 H_1[k]+2 H_2[k]+2 \beta  H_3[k]+2 \beta  H_4[k]+H_5[k]. 
\end{equation}
Then we investigate $H_1[k]$,$H_2[k]$, $H_3[k]$, $H_4[k]$, and $H_5[k]$.

Due to the  $L$-Lipschitz continuity of $\boldsymbol{F}$, we have 
$$
\begin{aligned}
	& 2 \beta(\boldsymbol{z}[k+1]-\tilde{\boldsymbol{z}}[k])^\mathrm{T}(\boldsymbol{F}(\boldsymbol{z}[k])-\boldsymbol{F}(\tilde{\boldsymbol{z}}[k])) \\
	& \leq 2 \beta\|\boldsymbol{z}[k+1]-\tilde{\boldsymbol{z}}[k]\| \cdot\|\boldsymbol{F}(\boldsymbol{z}[k])-\boldsymbol{F}(\tilde{\boldsymbol{z}}[k])\| \\
	& \leq 2 \beta \cdot L\|\boldsymbol{z}[k+1]-\tilde{\boldsymbol{z}}[k]\| \cdot\|\boldsymbol{z}[k]-\tilde{\boldsymbol{z}}[k]\| \\
	& \leq\|\boldsymbol{z}[k+1]-\tilde{\boldsymbol{z}}[k]\|^2+\beta^2 L^2\|\boldsymbol{z}[k]-\tilde{\boldsymbol{z}}[k]\|^2 .
\end{aligned}
$$
Therefore, 
$
H_1[k] \leq-\left(1-\beta^2 L^2\right)\|\boldsymbol{z}[k]-\tilde{\boldsymbol{z}}[k]\|^2.
$

On the other hand, following the basic property of projection in (\ref{projection}), we have 
$$
\begin{aligned}
		& H_2[k] \leq 0,\quad H_5[k] \leq 0.
\end{aligned}
$$
Moreover,  because $\boldsymbol{z}^{\dag}\in \Upsilon $ is the fixed point of algorithm (\ref{dis-compact}),
we have $ H_3[k] \leq 0$. 
Also,   $ H_4[k] \leq 0$ since $F(\boldsymbol{z})$ is monotone. 

Therefore, $\Delta[k]\leq0$, which  indicates that $\boldsymbol{z}[\cdot]$ is bounded.  Thus, there exists a subsequence
$k_{p}>0$, $k_{p} \rightarrow \infty$ as ${p} \rightarrow \infty$ and a point $\boldsymbol{z}^{\Diamond} \in  \boldsymbol{\varXi}$ such that
$$
\lim _{p \rightarrow \infty} \boldsymbol{z}\left[k_p\right]=\boldsymbol{z}^{\Diamond}.
$$ 
Then we prove that $\boldsymbol{z}^{\Diamond} $ is a fixed point of (\ref{dis-compact}), i.e, $\boldsymbol{z}^{\Diamond}\in \Upsilon $.

Since $ \Delta[k]\leq -(1-\beta^2L^2)\|\boldsymbol{z}[k]-\tilde{\boldsymbol{z}}[k]\|^2$,
$$
0 \leq \sum\nolimits_{k=1}^{\infty}\left(1-\beta^2 L^2\right)\|\boldsymbol{z}[k]-\tilde{\boldsymbol{z}}[k]\|^2 \leq \Delta[0],
$$
which implies 
$
\lim \nolimits_{p \rightarrow \infty} \|\boldsymbol{z}[k]-\tilde{\boldsymbol{z}}[k] \|=0.
$
As a result, 
$$
\lim _{p \rightarrow \infty} \tilde{\boldsymbol{z}}[k_p]=\boldsymbol{z}^{\Diamond}. \vspace{-0.2cm}
$$
Moreover,
it follows from $\tilde{\boldsymbol{z}}[k]=\Pi_{{\boldsymbol{\varXi}}}[\boldsymbol{z}[k]-\beta \boldsymbol{F}(\boldsymbol{z}[k])]$ that
$
\boldsymbol{z}^{\Diamond}=\Pi_{{\boldsymbol{\varXi}}}(\boldsymbol{z}^{\Diamond}-\beta \boldsymbol{F}(\boldsymbol{z}^{\Diamond})).
$
According to (\ref{fix2}), $\boldsymbol{z}^{\Diamond}\in \Upsilon$.

Finally, we show that $\boldsymbol{z}^{\Diamond}$ is not only a cluster point but also 
the unique limit point of $\boldsymbol{z}[k]$ as $k \rightarrow \infty$. Define $\bar{\Delta}[k]=\|\boldsymbol{z}[k+1]-\boldsymbol{z}^{\Diamond}\|^2-\|\boldsymbol{z}[k]-\boldsymbol{z}^{\Diamond}\|^2$.  With similar analysis, we prove $\bar{\Delta}[k]\leq0$. In other words, $\|\boldsymbol{z}[k]-\boldsymbol{z}^{\Diamond}\|^2$ is monotonically decreasing with 
\begin{equation}
	\lim _{k \rightarrow \infty} \|\boldsymbol{z}[k]-\boldsymbol{z}^{\Diamond}\|^2=\lim _{p \rightarrow \infty} \|\boldsymbol{z}[k_p]-\boldsymbol{z}^{\Diamond}\|^2=0.
\end{equation}
which yields   (\ref{lim}). Moreover, 
if the dual relation \eqref{dual-rela} is satisfied, then the convergent point $\boldsymbol{x}^{\Diamond}$ is indeed the global NE. 
This completes the proof. \hfill $\square$
\vspace{-0.25cm}


\section{Numerical Experiments}\label{sim}

In this section, 
we examine the effectiveness of Algorithm 1 and Algorithm 2 for seeking
 NE of the SNL problem.  

\begin{table*}
	\centering
	\footnotesize
	\caption{MLE of four methods  with different numbers of non-anchor nodes and anchor nodes.}
		\setlength\tabcolsep{15pt}
		\renewcommand\arraystretch{1.3}
		\begin{tabular}{c|c|c|c|c}
			\hline
			\hline
			\specialrule{0em}{0.2pt}{0.5pt}
 Node	numbers &$ M=10,N=10 $ &$ M=18,N=30 $  &$ M=30,N=70 $& $ M=40,N=100 $  \\ 
   \hline	
			Algorithm 2 &0.0213 &0.0164  &0.0153& 0.0147 \\ \hline
   Alg-RBR \cite{jia2013distributed} &0.0508 &0.0915  &0.1022& 0.684  \\ \hline
			Alg-SDP \cite{lui2008semi}&0.0773&0.4835  & 0.9729& 1.664 \\  \hline
			Alg-ARMA \cite{wan2019sensor}&0.0681 &0.6124 &0.8473 &1.254  \\\hline
			\hline
	\end{tabular}
	\label{tab61}
\end{table*}




We first  consider a two-dimensional case \cite{calafiore2010distributed} with $N = 7$ non-anchor nodes and $M=3$ anchor nodes located on the unit squares $[0,1]\times [0,1]$ according to the configuration
shown in Fig. \ref{fig4}(a). 
The anchor nodes  are shown by red asterisks, which are located on  the external vertices of the formation. The
true locations of non-anchor nodes are shown by  
blue diamonds.  The dashed lines indicate that nodes are connected. The initial distribution of non-anchor nodes is randomly generated, shown by purple asterisks.  The  localization effectiveness is evaluated by
the \textit{mean localization error}:
\begin{equation*}
  \operatorname{MLE}=\frac{1}{N} \sqrt{\sum\nolimits_{i=1}^{N}\|x_{i}-x_{i}^{\Diamond}\|^2}.
\end{equation*}
 We employ Algorithm 1 to solve this problem. Assign $\alpha_k=\frac{0.0637}{\sqrt{k}}$ as the step size. Also,  
set the  tolerance $t_{tol}= 10^{-3}$  and the terminal criterion 
$ \|\boldsymbol{x}[k+1]-\boldsymbol{x}[k]\| \leq t_{tol},\; \|{\sigma}[k+1]-{\sigma}[k+1]\| \leq t_{tol}.  $
	\begin{figure}[ht]
			\hspace{-1.7cm}
			\centering	
			\subfigure[Initial distribution]{
				\begin{minipage}[t]{0.45\linewidth}
					\centering
					\includegraphics[width=4.6cm]{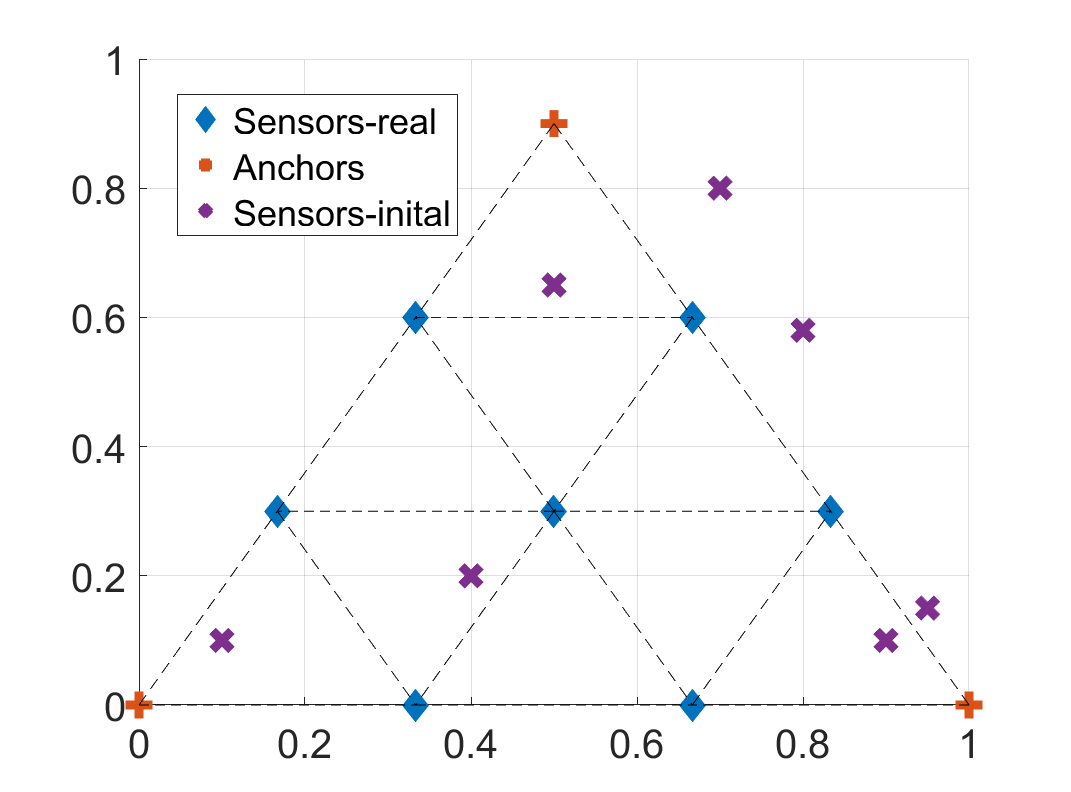}
				\end{minipage}%
			}%
			\hspace{0.1cm}
			\subfigure[Final localization results]{
				\begin{minipage}[t]{0.45\linewidth}
					\centering
					\includegraphics[width=4.6cm]{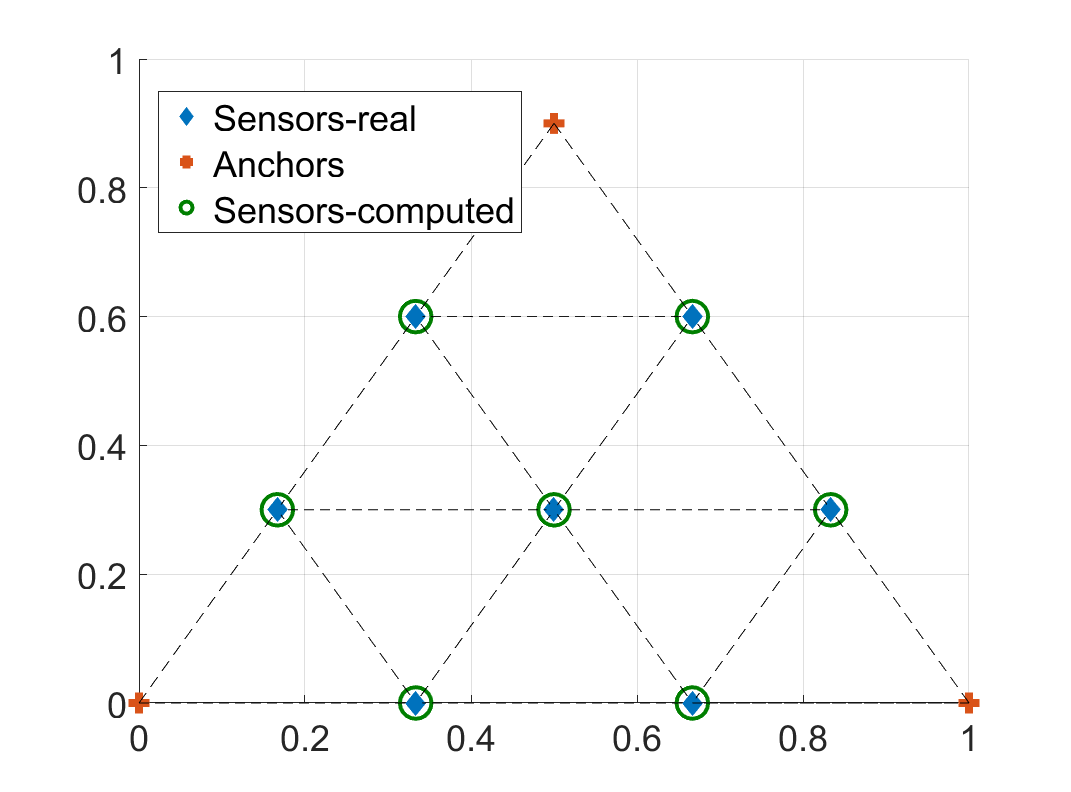}
				\end{minipage}%
			}%
			\centering
			\caption{The initial distribution and final localization results.} 
   \vspace{-0.1cm}
			\label{fig4}
		\end{figure}

Fig. \ref{fig4}(b) shows the final localization by Algorithm 1, where the anchor nodes and
the true locations of non-anchor nodes are also shown by red
asterisks and blue diamonds, and the computed locations are shown
by green circles. As shown in Fig. \ref{fig4}(b),  the computed
location results match the true locations. This  indicates that
our approach ensures the location  accuracy of non-anchor nodes and is
able to find the global NE under this non-convex scenario. 
%
Also,  Fig. \ref{fig12}  shows the trajectories of all non-anchor nodes’ strategies in Algorithm 1 with respect to one certain dimension. 
Since the set $ \mathbb{E}^{+}$ is nonempty and the identification and verification condition in Theorem \ref{cent} is satisfied, 
all non-anchor nodes can find their appropriate localization on account of the convergence, 
which actually serves as the desired global NE.

  \begin{figure}[ht]
	\centering	
	\includegraphics[scale=0.2]{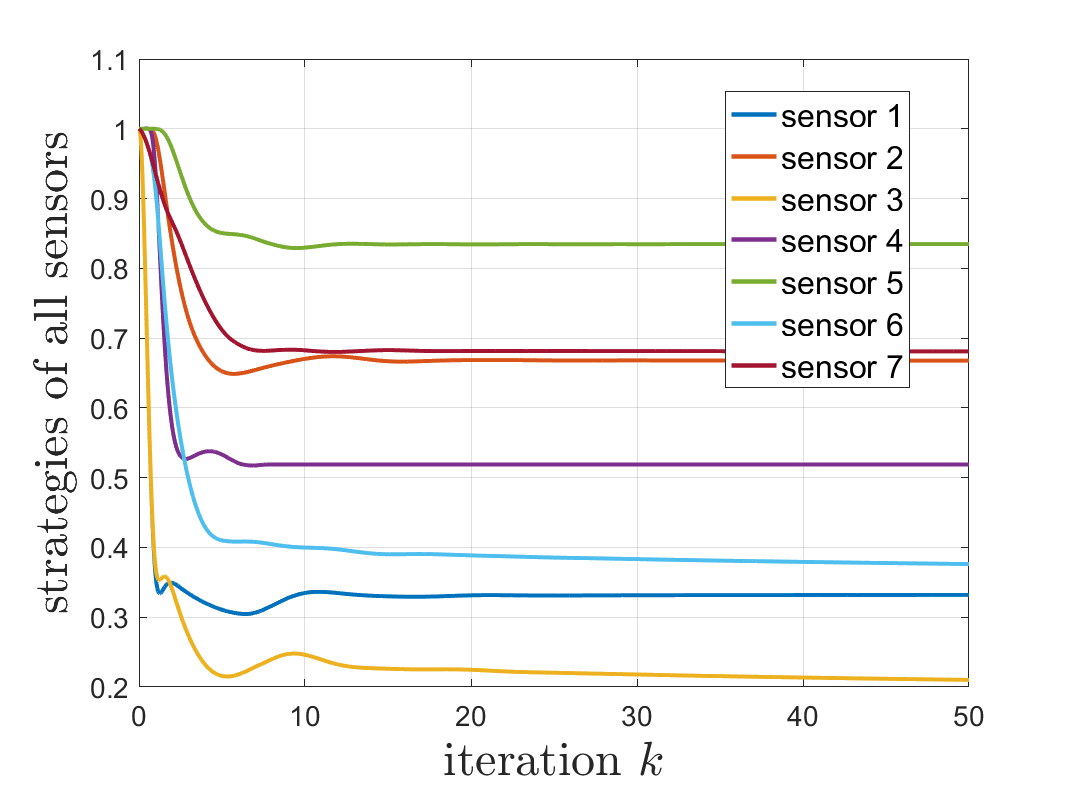}\\
 \vspace{-0.2cm}
	\caption{Convergence of all non-anchor nodes' decisions}
	\label{fig12}
	\vspace{-0.3cm}
\end{figure}
Moreover, we show the effectiveness of Algorithm 1 for SNL problems with different node configurations. The first one is with $N = 2$ non-anchor nodes and $M=4$ anchor nodes, which are randomly generated on the unit squares $[-2,2]\times [-2,2]$. The second one is with $N = 50$ non-anchor nodes and $M=18$ anchor nodes, which are randomly generated on the unit squares $[-3,3]\times [-3,3]$.  Fig. \ref{figt4} presents the computed sensor location results in these cases. Clearly,  Algorithm 1  can localize all sensors in either small or large sensor networks. 


  \begin{figure}[ht]
			\centering	
			\subfigure[$N\!=\!2$, $M\!=\!4$]{
				\begin{minipage}[t]{0.45\linewidth}
					\centering
					\includegraphics[width=4.5cm]{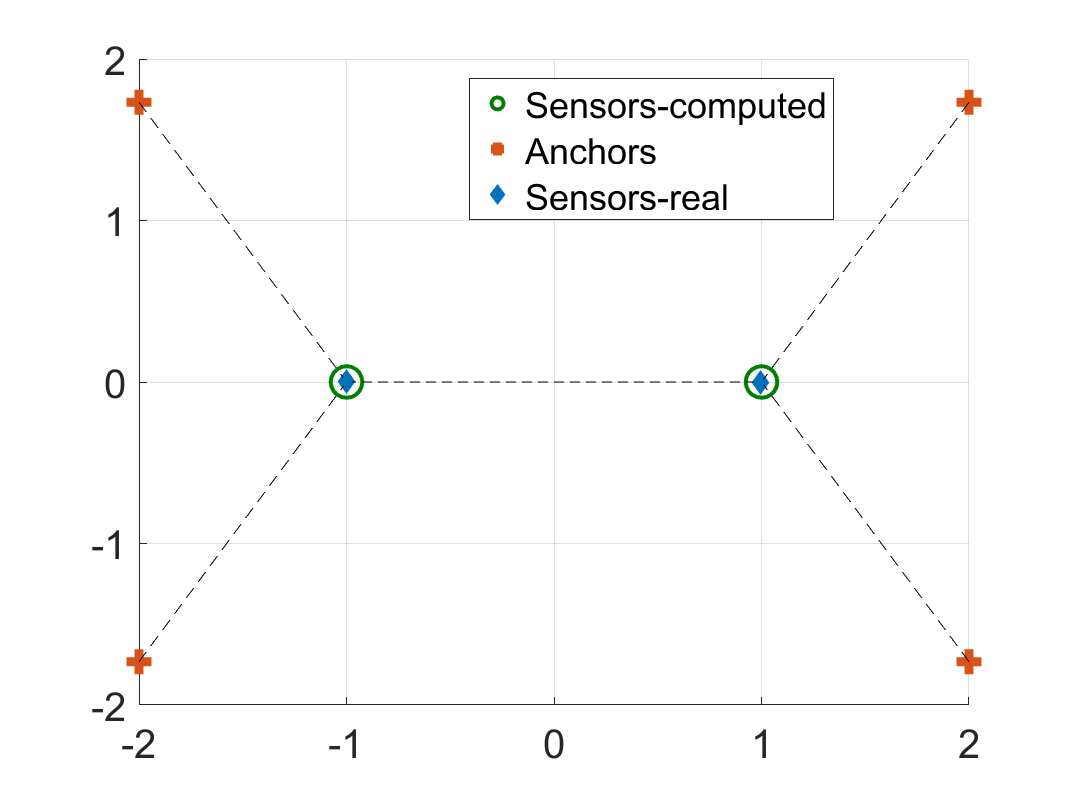}
     \vspace{-0.2cm}
				\end{minipage}%
			}%
			\hspace{0.3cm}
			\subfigure[$N\!=\!50$, $M\!=\!18$]{
				\begin{minipage}[t]{0.45\linewidth}
					\centering
					\includegraphics[width=3.9cm]{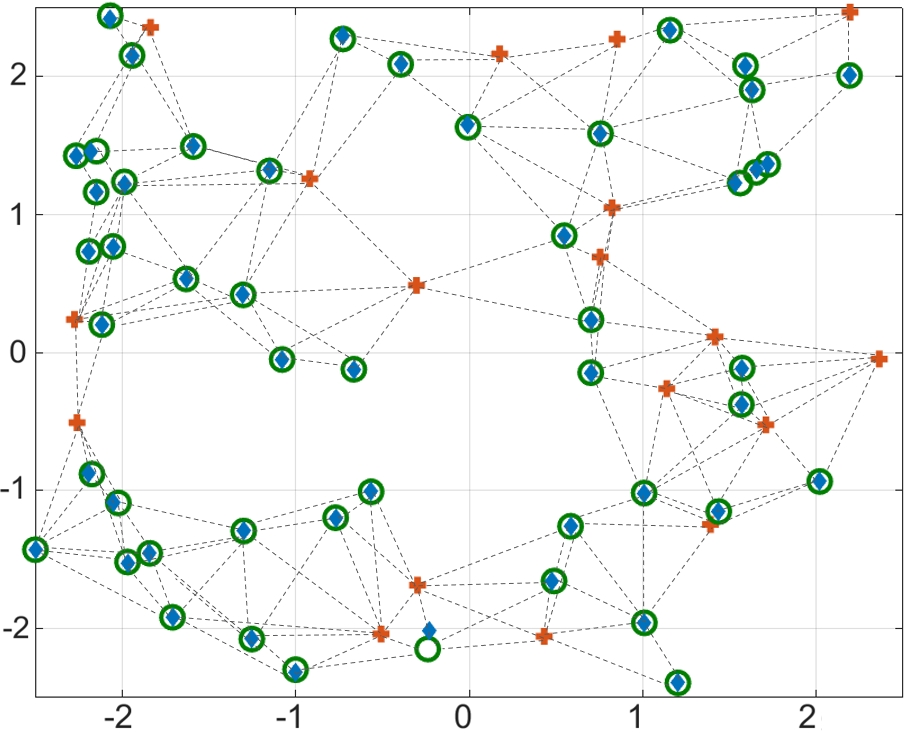}
				\end{minipage}%
			}%
			\centering
   \vspace{-0.15cm}
			\caption{Computed sensor location results with different configurations.} 
			\label{figt4}
   \vspace{-0.4cm}
		\end{figure}

On the other hand, 
we compare the performance of Algorithm 2 with those of several distributed algorithms for solving the SNL problem under both noise-free and noisy situations.
These compared algorithms include distributed SDP-based algorithm \cite{lui2008semi}, distributed ARMA-based algorithm \cite{wan2019sensor}, and distributed random best response (RBR) algorithm \cite{jia2013distributed}. We take the 
the UJIIndoorLoc data set. The UJIIndoorLoc data set was introduced in 2014 at the International Conference on Indoor Positioning and Indoor Navigation, which estimates user location based on building and floor. 
The UJIIndoorLoc data set is available on the UC Irvine Machine Learning Repository website \cite{23}.
For the  noise-free environment, we randomly generate four different node configurations 
on the unit square $[-5,5]\times [-5,5]$ and record the value of MLE under different algorithms in Table I.  Table I reflects that  Algorithm 2 achieves a higher localization accuracy than other methods. Moreover,  with the expansion of network size, only Algorithm 2 still maintains MLE in a tolerance error
range, while other methods can not guarantee this.  Specifically, Fig. \ref{figg4} shows the computed results  from four methods for case $N=30$, $M=18$. In Fig.  \ref{figg4}, Algorithm 2   localizes all sensors, while other methods  are accompanied by  deviation errors from the true locations.

    \begin{figure}[ht]
			\centering	
			\subfigure[results from Alg. 2]{
				\begin{minipage}[t]{0.49\linewidth}
					\centering
					\includegraphics[width=3.7cm]{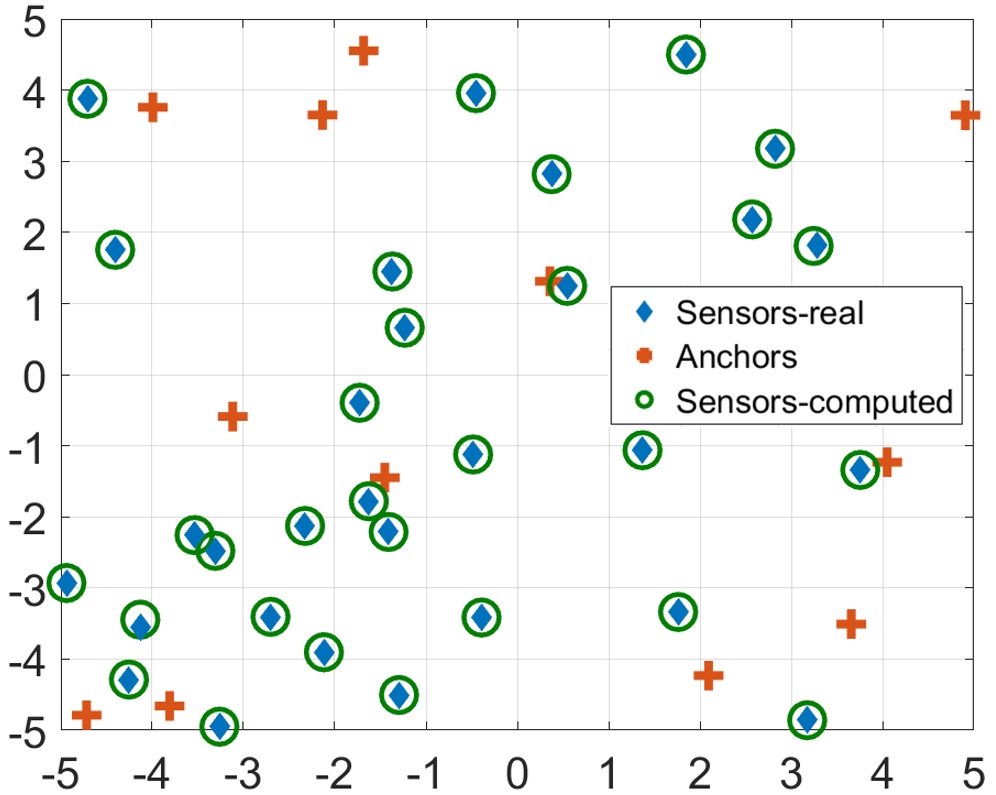}
				\end{minipage}%
			}%
			\subfigure[results from Alg-RBR  \cite{jia2013distributed}]{
				\begin{minipage}[t]{0.49\linewidth}
					\centering
					\includegraphics[width=3.7cm]{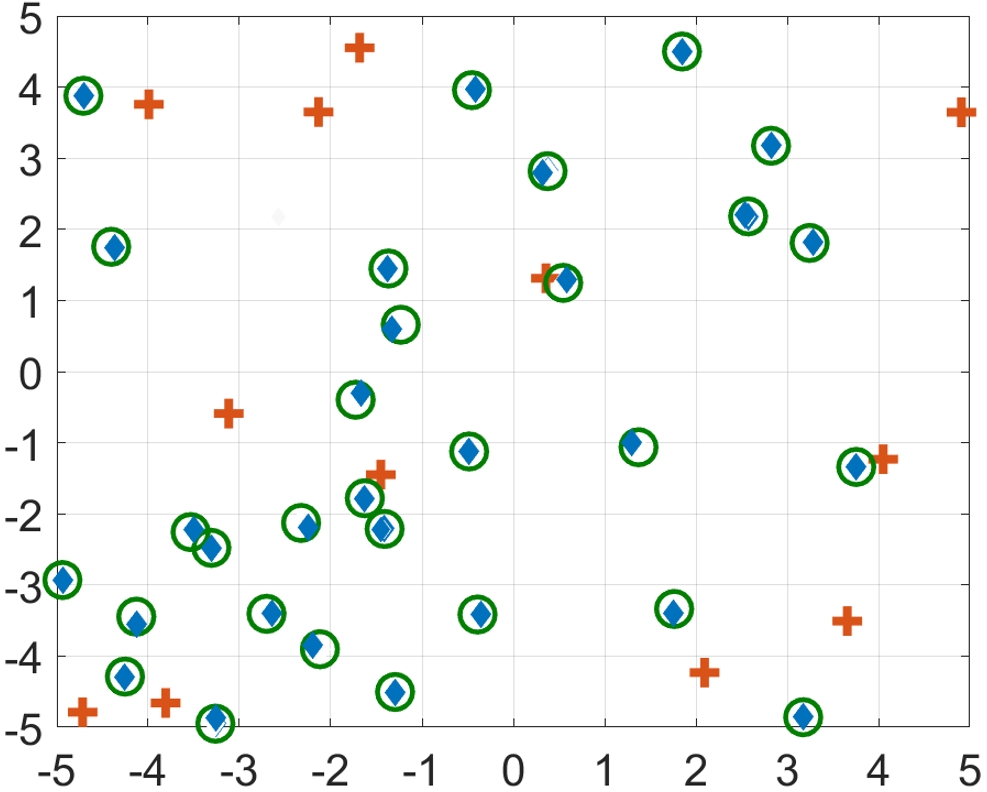}
				\end{minipage}%
			}%
   \\
   \subfigure[results from Alg-SDP \cite{lui2008semi}]{
				\begin{minipage}[t]{0.49\linewidth}
					\centering
					\includegraphics[width=3.7cm]{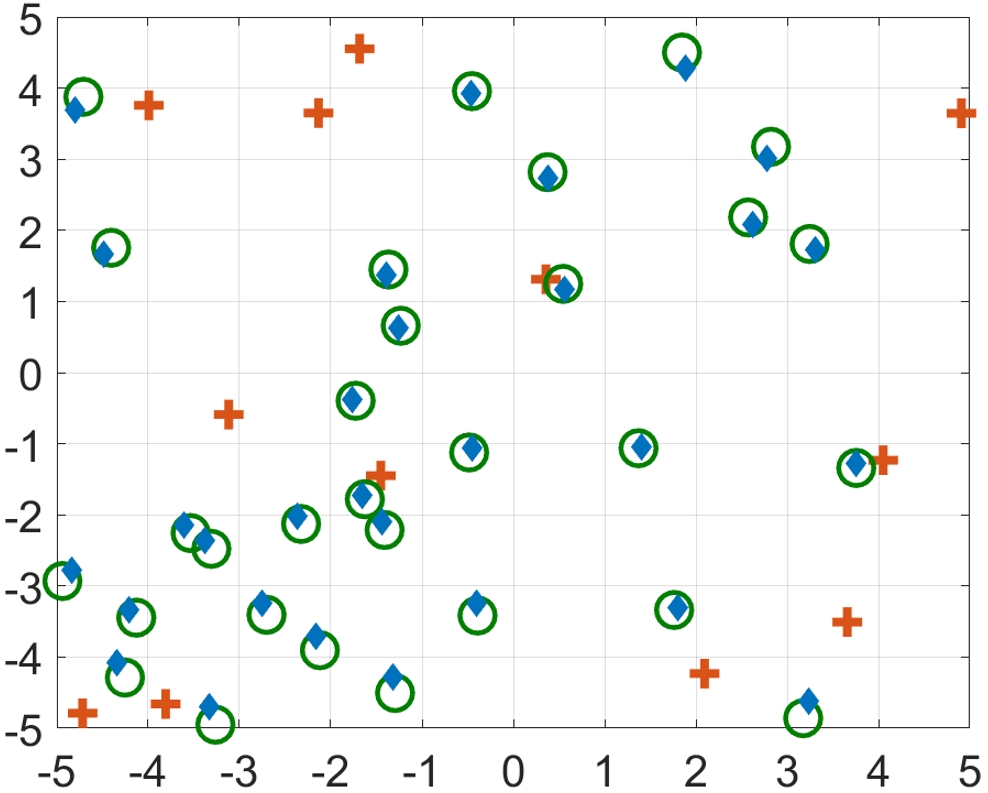}
     \vspace{-0.2cm}
				\end{minipage}%
			}%
			\subfigure[results from Alg-ARMA \cite{wan2019sensor}]{
				\begin{minipage}[t]{0.49\linewidth}
					\centering
					\includegraphics[width=3.7cm]{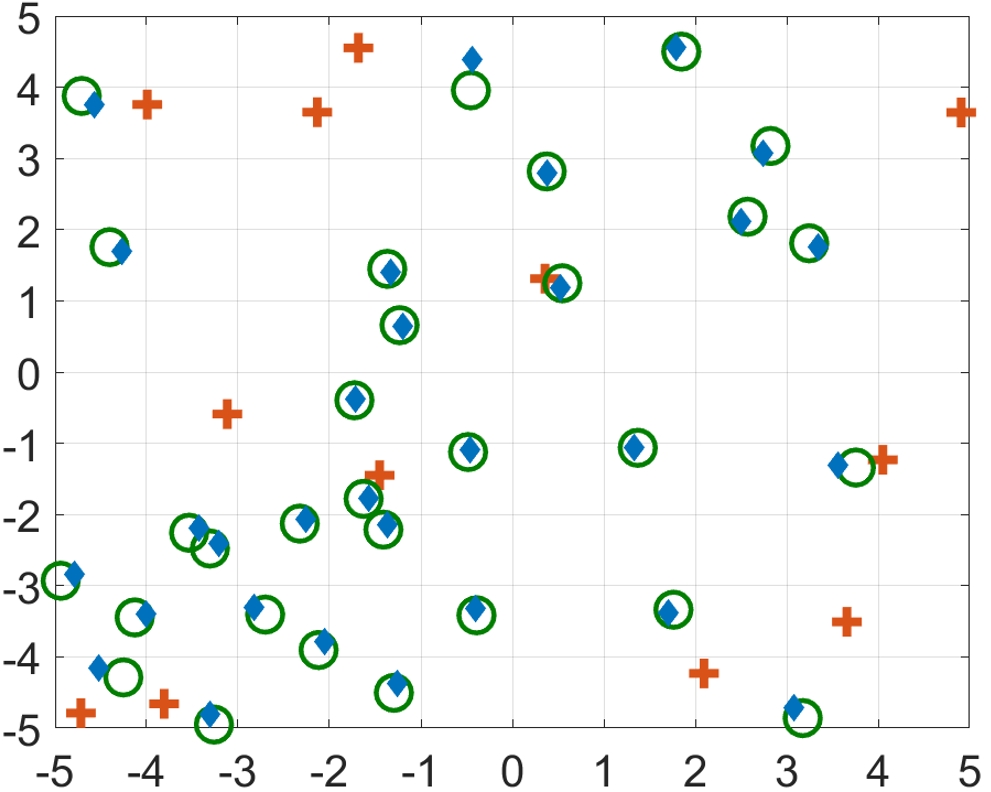}
				\end{minipage}%
			}
			\centering
			\caption{Computed sensor location results with different configurations.} 
			\label{figg4}
    \vspace{-0.3cm}
		\end{figure}









Even if in this paper we assume that no noise affects the measurement, for completeness of the simulation examples, we consider a case with noisy measurements. 
Let  the range measurement be $d_{ij}=\bar{d}_{ij}*(1+N(0,0.001))$ and $e_{il}=\bar{d}_{ij}*(1+N(0,0.001))$, where $\bar{d}_{ij}$ and $\bar{e}_{il}$ are the actual distances between the two nodes, and $N(0,0.001)$  is a random variable. 
A network of $N=18$ non-anchor nodes and $M=6$ anchor nodes from the UJIIndoorLoc data set are generated in the square area $[-2.5,2.5]\times [-2.5,2.5]$.
 The numerical results obtained by the four methods 
are plotted in Fig. \ref{fige4}, where Algorithm 2 localizes all these sensors with much higher accuracy than other methods.  Also, this indicates that our approach can be extended to noisy cases.

   \begin{figure}[ht]
			\centering	
			\subfigure[results from Alg. 2]{
				\begin{minipage}[t]{0.49\linewidth}
					\centering
					\includegraphics[width=3.8cm]{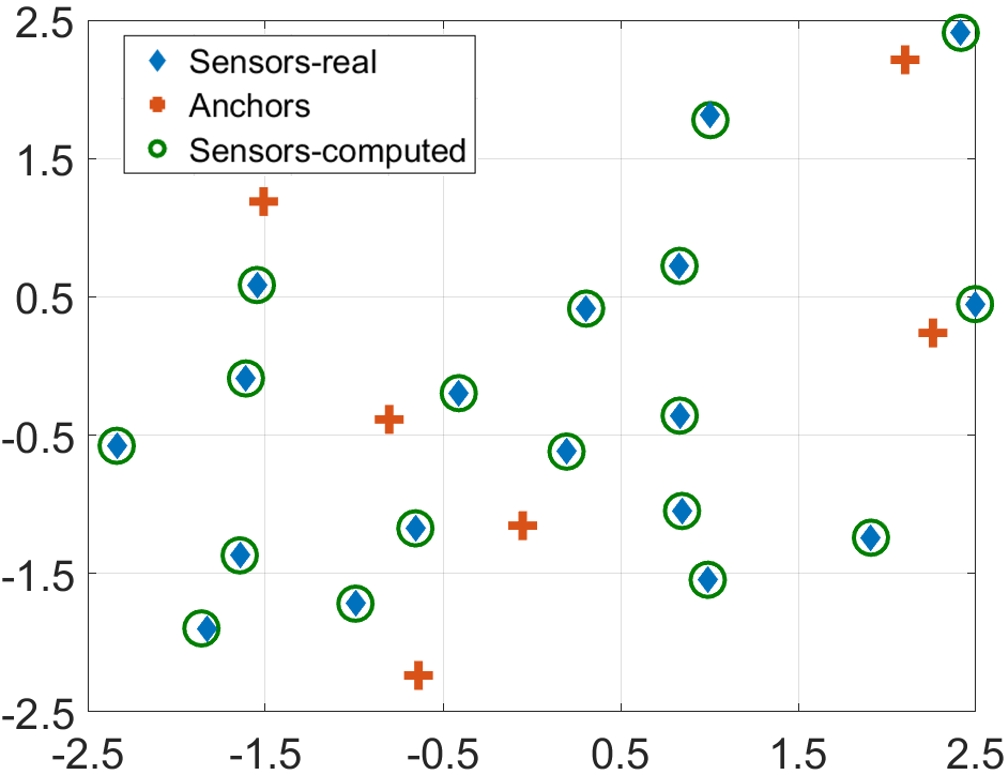}
     \vspace{-0.2cm}
				\end{minipage}%
			}%
			\subfigure[results from Alg-RBR \cite{jia2013distributed}]{
				\begin{minipage}[t]{0.49\linewidth}
					\centering
					\includegraphics[width=3.8cm]{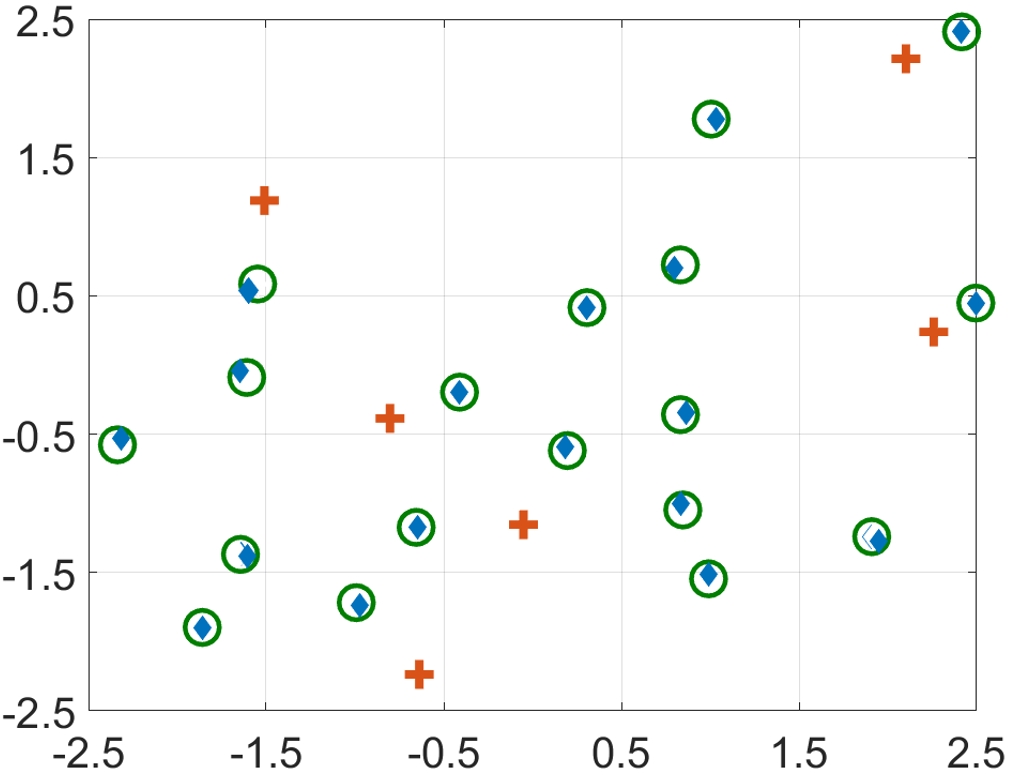}
				\end{minipage}%
			}%
   \\
   \subfigure[results from Alg-SDP \cite{lui2008semi}]{
				\begin{minipage}[t]{0.49\linewidth}
					\centering
					\includegraphics[width=3.8cm]{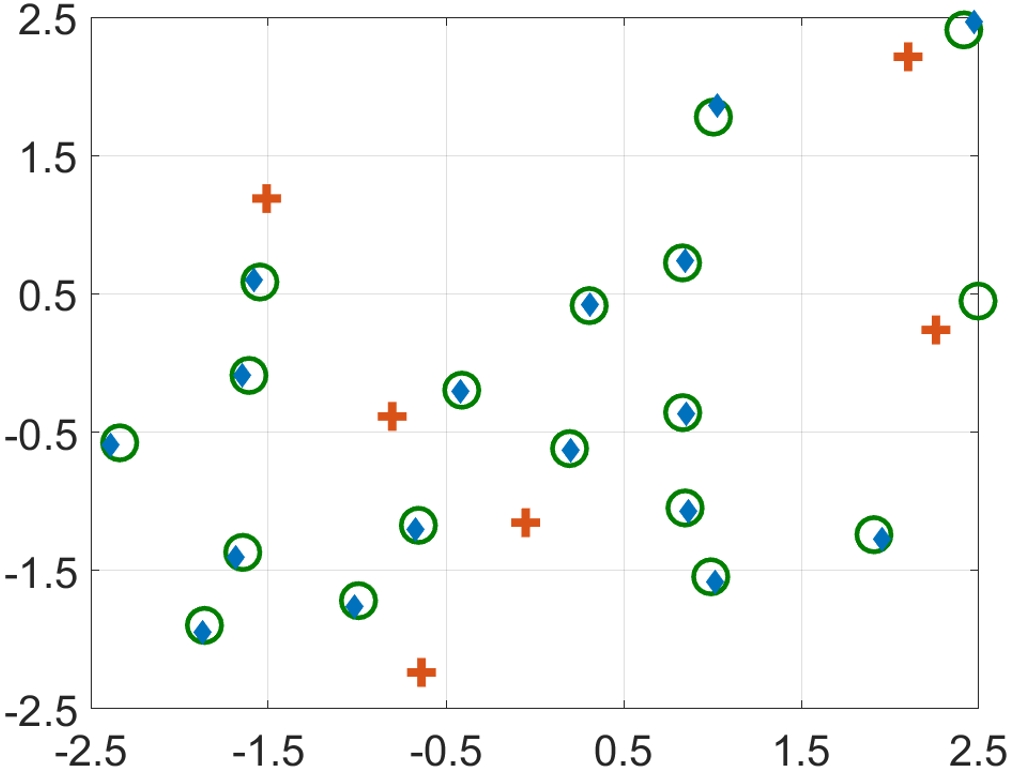}
     \vspace{-0.2cm}
				\end{minipage}%
			}%
			\subfigure[results from Alg-ARMA \cite{wan2019sensor}]{
				\begin{minipage}[t]{0.49\linewidth}
					\centering
					\includegraphics[width=3.8cm]{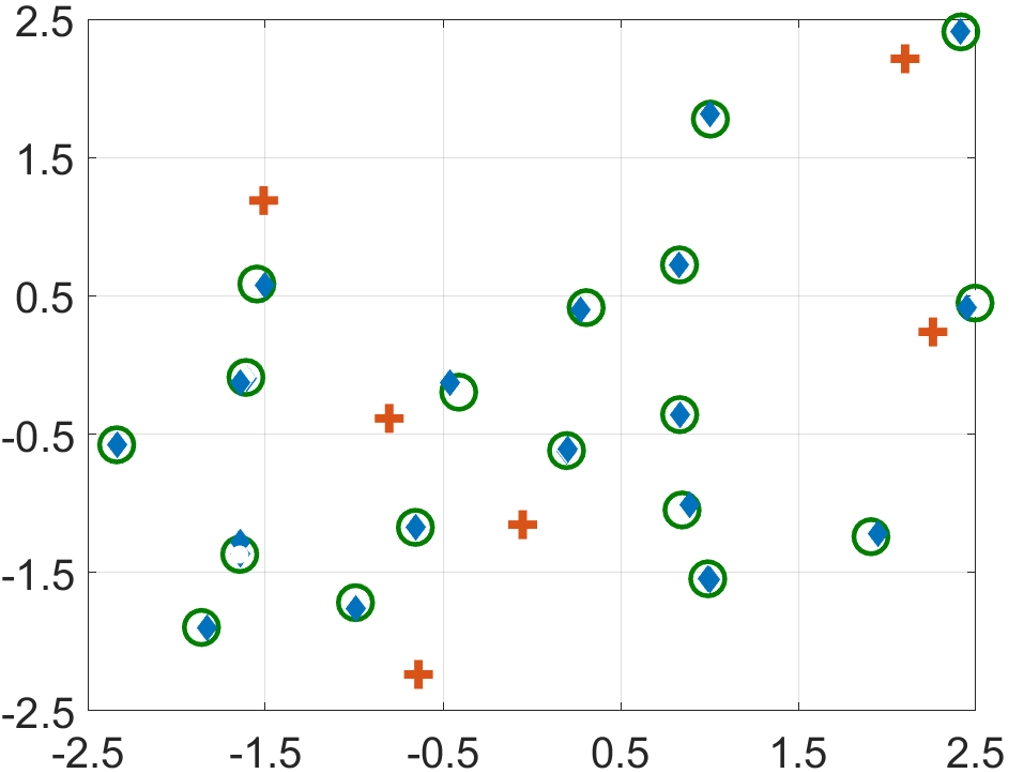}
				\end{minipage}%
			}
			\centering
			\caption{Computed sensor location results with different configurations.} 
			\label{fige4}
   \vspace{-0.4cm}
		\end{figure}

\section{Concluding Remarks}
\vspace{-0.1cm}
In this paper, we have focused on  the non-convex SNL problem. We have presented novel results on both the identification and verification condition of the global solution and the position-seeking algorithms. By  formulating the SNL problem as a potential game, we have shown that the  NE  is equivalent to the global network's solution. Then, based on the canonical duality theory, we have proposed  a primal-dual algorithm  to compute the stationary point of a complementary dual problem  with an $\mathcal{O}(1/\sqrt{k})$ rate, which  actually induces the 
NE  if a duality relation can be checked. Moreover,  we have assigned   the algorithm to distributed implementation by virtue of sliding mode control and extra-gradient methods,  and proposed another distributed approach   to find the global NE, followed with  the analysis of its global convergence. Finally, the computational efficiency of our algorithm has been illustrated by several experiments. 

Future research efforts will focus on  extending our current results to more complex scenarios  such as i) generalizing the model to cases with measurement noise,  ii) quantitatively analyzing the influence of the uncertainty  on the equilibrium and the robustness of the algorithm, and iii) looking for weaker graph conditions.

 \begin{appendices}
\section{Proof of Lemma \ref{r11}}\label{r122}
\vspace{-0.13cm}
Recall the formulation of the non-convex function $\Phi$. If all the non-anchor nodes are localized accurately such that $\boldsymbol{x}^{\star}$ satisfies $\|x_{i}^{\star}-x_{j}^{\star}\|^2-d_{i j}^2=0$ and $\|x_{i}^{\star}-a_{l}\|^2-e_{i l}^2=0$ for any  $(i,j)\in \mathcal{E}_{ss}$ and $(i,l)\in \mathcal{E}_{as}$, then 
$\Phi(x_{1}^{\star},\!\dots,\! x_{N}^{\star})$ tends to  zero. Thus  $\boldsymbol{x}^{\star}$ is the global minimum of $\Phi$. 
Moreover, referring to  \cite{eren2004rigidity, aspnes2006theory}, Assumption 1 guarantees that the global minimum of $\Phi$ is unique and corresponds to the actual position profile $\boldsymbol{x}^{\star}$. This completes the proof.  \hfill $\square$

\section{Proof of Proposition \ref{t3}}\label{t122}

Given that the non-anchor nodes’ payoff functions $J_{i}$ are continuously differentiable, it follows from \cite{monderer1996potential} that the condition (\ref{pp1}) in Definition \ref{d2} is  equivalent to 
$$
\frac{\partial J_{i}(x_{i}, \boldsymbol{x}_{-i})}{\partial x_i}=\frac{\partial \Phi(x_{i}, \boldsymbol{x}_{-i})}{\partial x_i}, \quad \text { for } i\in\{1,\dots,N\}.
$$  
By computing the partial derivatives of $J_{i}$ and $\Phi$, respectively, 
we can verify that $G$ is a potential game.
 
Moreover, 
referring to \cite[Corollary 2.1, Theorem 2.3]{la2016potential}, we can derive that 
in a potential game,   a global NE corresponds exactly to a global
minimum of the potential function. 
By Lemma \ref{r11}, we yield that $\boldsymbol{x}^{\Diamond}$ is unique and is equal to $\boldsymbol{x}^{\star}$.
\hfill $\square$

\section{Details of domain $\mathbb{E}^{+}$}\label{aa}

		Define the second-order partial derivative of  $\Gamma(  \boldsymbol{x},{\sigma})$ in $\boldsymbol{x}$ as 
		$
			P(\sigma)= \nabla^{2}_{\boldsymbol{x}}\Gamma=
			2\left(\operatorname{diag}\left(V_1\left(\sigma^{s}\right)\right)+\operatorname{diag}\left(V_2\left(\sigma^{e}\right)\right)+G\left(\sigma^s\right)\right)
		$
		with $ V_1\left(\sigma^{s}\right),V_2\left(\sigma^{e}\right)\in \mathbb{R}^{Nn}$ and $G\left(\sigma^s\right)\in \mathbb{R}^{Nn\times Nn}$, where
		$$
		V_1\left(\sigma^{s}\right)=\left[\begin{array}{c}
			\sum_{j=1}^N \sigma_{1 j}^s \!\\
			\vdots \\
			\sum_{j=1}^N \sigma_{1 j}^s \\
			\vdots \\
			\sum_{j=1}^N \sigma_{N j}^s \\
			\vdots \\
			\sum_{j=1}^N \sigma_{N j}^s
		\end{array}\right], V_2\left(\sigma^{e}\right)=\left[\begin{array}{c}
			\sum_{k=1}^{M} \sigma_{1 k}^a \\
			\vdots \\
			\sum_{k=1}^{M} \sigma_{1 k}^a \\
			\vdots \\
			\sum_{k=1}^{M} \sigma_{N k}^a \\
			\vdots \\
			\sum_{k=1}^{M} \sigma_{N k}^a
		\end{array}\right] \text {, }
		$$
		$$
		G\left(\sigma^s\right)=\left[\begin{array}{ccc}
			-\sigma_{11}^s I_{n} & \dots & -\sigma_{1 N}^s I_{n} \\
			\vdots & \vdots & \vdots \\
			-\sigma_{N 1}^s I_{n} & \dots & -\sigma_{N N}^s I_{n}
		\end{array}\right],
		$$
		where  
		$\sigma_{i j}^s=\sigma_{ji}^s$, $\sigma_{i j}^s=0$ if $(i,j)\notin\mathcal{E}_{ss}$ and $\sigma_{i l}^a=0$ if $(i,l)\notin\mathcal{E}_{as}$.
		{
			It is obvious that $\nabla^{2}_{\boldsymbol{x}}\Gamma$ is   both $x_i$-free and $\boldsymbol{x_{-i}}$-free,  and can be regarded as a linear combination for the elements of $\sigma$. } 
Moreover,  we can verify that 
$\boldsymbol{0}_{q}\in\{\sigma : 
	P(\sigma) \succeq 0
	\}$ by taking $\sigma_{i j}^s=0 $ and $\sigma_{i l}^a=0$ for all $(i,j)\in\mathcal{E}_{ss}$ and $(i,l)\notin\mathcal{E}_{as}$. Also, 
 it is easy to find
 $\boldsymbol{0}_{q}\in\Theta^{*}$. Thus, $\mathbb{E}^{+}= \Theta^{*}\cap \{\sigma : 
	P(\sigma) \succeq 0
	\}$ is non-empty. 
For any given $\sigma^{'}, \sigma^{''}\in \mathbb{E}^{+}$, we have 
$
\omega P(\sigma^{'}) \succeq 0, (1-\omega) P(\sigma^{''}) \succeq 0,\, \forall \omega\in [0,1].
$
Thus, for $\omega\in [0,1]$,
$$
\omega P(\sigma^{'}) +(1-\omega) P(\sigma^{''}) = P(\omega\sigma^{'}+(1-\omega)\sigma^{''} ) \succeq 0.
$$
   This implies that $\mathbb{E}^{+}$ is convex. \hfill $\square$

\vspace{-0.2cm}
\section{Proof of Theorem  \ref{t2}}\label{bb}

For  convenience, we define the following compact profiles for variables in \eqref{ci}, $\Lambda\left(\boldsymbol{x}\right)=\operatorname{col}\left\{
	\psi^{s}(\boldsymbol{x}), 	\varphi^{a}(\boldsymbol{x})
	\right\}$. We have $\xi=\operatorname{col}\{\xi^{s},\xi^{e} \}=\Lambda\left(\boldsymbol{x}\right)$, and $\xi\in\Theta \subseteq \mathbb{R}^{q}$,
where $\Theta=\Theta_s\times\Theta_a$
and $q=\left|\mathcal{E}_{ss}\right|+\left|\mathcal{E}_{as}\right|$. Following the definitions of \eqref{cl} and \eqref{cr}, we denote
\begin{equation*}
	\begin{aligned}
		&\Psi\left(\xi\right)= \Psi_{s}({\xi}^s)+\Psi_{a}({\xi}^a),\;\Psi^{*}\left(\sigma\right)= \Psi_{s}^*({\sigma}^s)+\Psi_{a}^*({\sigma}^a).
	\end{aligned}
\end{equation*} 
Then it follows from the canonical duality relation  \eqref{dual} that
$$
\sigma=\nabla\Psi({\xi})\Leftrightarrow \xi=\nabla\Psi^*({\sigma})\Leftrightarrow 	\xi^{\operatorname{T}} \sigma=\Psi(\xi) +\Psi^*(\sigma). 
$$
Thus, the  
function $\Gamma$ in \eqref{complementary} can be expressed as
$$
\Gamma\left(\boldsymbol{x},\sigma\right)	=\sigma^{\operatorname{T}} \Lambda\left(\boldsymbol{x}\right)-\Psi^{*}\left(\sigma\right).\vspace{-0.2cm}
$$
		For the conjugate gradient of 
  $\Psi$,
  denote  
		\begin{equation*}
			\nabla\Psi^{*}\left(\sigma\right)=
			\left[\begin{array}{c}
				\operatorname{col}\{\frac{1}{2}(\sigma_{ij}^{s})+d_{ij}^2 \}_{(i,j)\in\mathcal{E}_{ss}}\\ \operatorname{col}\{\frac{1}{2}(\sigma_{il}^{a})+e_{il}^2\}_{(i,l)\in\mathcal{E}_{as}}
			\end{array}\right]
		\end{equation*} 
		Then
		the partial derivative profile  of $\Gamma$ in $\boldsymbol{x}$ is expressed as   $\nabla_{\boldsymbol{x}} \Gamma\left( \boldsymbol{x},\sigma\right) =\nabla_{\boldsymbol{x}}\,\sigma^{\operatorname{T}} \Lambda\left(\boldsymbol{x}\right).$
	Also,
		the partial derivative profile of $\Gamma$ in $\sigma$ is expressed as $	\nabla_{\sigma} \Gamma\left( \boldsymbol{x},\sigma\right)  =
			\Lambda\left(\boldsymbol{x}\right)-\nabla \Psi^{*}\left(\sigma\right).$

  
Now we prove  Theorem  \ref{t2}.
We first prove the sufficiency.  Under Assumption 1, if there exists ${\sigma}^{\Diamond}\in \mathbb{E}^{+}$ such that $(\boldsymbol{x}^{\Diamond},{\sigma}^{\Diamond}) \in \boldsymbol{\Omega}\times\mathbb{E}^{+} $ 
	is a stationary point of    $ \Gamma (\boldsymbol{x},\sigma)$, then 
\begin{subequations}\label{ddd}
\begin{align}
	&\mathbf{0}_{nN} \in   \nabla_{\boldsymbol{x}}\,\sigma^{\Diamond\operatorname{T}} \Lambda\left(\boldsymbol{x}^{\Diamond}\right)	
	+\mathcal{N}_{\boldsymbol{\Omega}}(\boldsymbol{x}^{\Diamond}),\\
	&\mathbf{0}_{q} \in -\Lambda(\boldsymbol{x}^{\Diamond})+\nabla \Psi^{*}(\sigma^{\Diamond})+ \mathcal{N}_{\mathbb{E}^{+}}(\sigma^{\Diamond}),\vspace{-0.2cm}
 \end{align}
\end{subequations} 
where $\mathcal{N}_{\boldsymbol{\Omega}}(\boldsymbol{x}^{\Diamond})$ is the normal cone at point $\boldsymbol{x}^{\Diamond}$ on set $\boldsymbol{\Omega}$, with a similar definition for the normal cone $\mathcal{N}_{\mathscr{E}^{+}}(\sigma^{\Diamond})$.
Moreover, if the canonical duality relation satisfies $\sigma_{i j}^{s\Diamond}=\nabla_{\xi_{i j}^s}\Psi_{s}({\xi}^s)|_{\xi_{i j}^s=\|x_{i}^{\Diamond}-x_{j}^{\Diamond}\|^2}, \forall (i, j) \in \mathcal{E}_{ss} $ and
	$\sigma_{i l}^{a\Diamond}=\nabla_{\xi_{i l}^a}\Psi_{a}({\xi}^a)|_{\xi_{i j}^a=\|x_{i}^{\Diamond}-a_{l}\|^2}, \forall (i, l) \in \mathcal{E}_{as}$, we have   
\begin{align*}
	\sigma^{\Diamond}=\nabla \Psi(
	\Lambda(\boldsymbol{x}^{\Diamond}))\Longleftrightarrow \;\Lambda(\boldsymbol{x}^{\Diamond})=\nabla \Psi^{*}(\sigma^{\Diamond}). 
\end{align*}
On this basis, (\ref{ddd}b) is transformed into
$
\Lambda(\boldsymbol{x}^{\Diamond})=\nabla \Psi^{*}(\sigma^{\Diamond}).
$
By substituting $\sigma^{\Diamond}$ with $\nabla \Psi(
	\Lambda(\boldsymbol{x}^{\Diamond}))$,
(\ref{ddd}a) becomes
\begin{equation}\label{chain}
    \mathbf{0}_{nN} \in   \nabla_{\boldsymbol{x}}\,\Psi(
	\Lambda(\boldsymbol{x}^{\Diamond}))^{\operatorname{T}} \Lambda\left(\boldsymbol{x}^{\Diamond}\right)	
	+\mathcal{N}_{\boldsymbol{\Omega}}(\boldsymbol{x}^{\Diamond}).
\end{equation}
According to the chain rule, 
 $
\nabla \Psi(
	\Lambda(\boldsymbol{x}^{\Diamond}))^{\operatorname{T}} \Lambda\left(\boldsymbol{x}^{\Diamond}\right)	= \nabla_{\boldsymbol{x}}\Phi(\boldsymbol{x}^{\Diamond}).
 $
 Therefore, \eqref{chain} is equivalent to 
\begin{equation}\label{ft41}
	 \mathbf{0}_{nN} \in   \nabla_{\boldsymbol{x}}\Phi(\boldsymbol{x}^{\Diamond})+\mathcal{N}_{\boldsymbol{\Omega}}(\boldsymbol{x}^{\Diamond}).\vspace{-0.13cm}
\end{equation}
According to the definition of potential game, \eqref{ft41} implies 
\begin{equation}\label{dde}
	\mathbf{0}_{n} \in \nabla_{x_{i}} J_{i}(x_{i}^{\Diamond}, \boldsymbol{x}_{-i}^{\Diamond})+\mathcal{N}_{\Omega_{i}}({x_{i}^{\Diamond}}).  \vspace{-0.13cm}
\end{equation}
 Moreover, when  $ \sigma\in\mathbb{E}^{+}  $,  the Hessian matrix satisfies
$
	\nabla_{\boldsymbol{x}}^{2} \Gamma(\boldsymbol{x},\sigma) \succeq \boldsymbol{0}_{nN},
 $
which indicates 
the convexity of $\Gamma(\boldsymbol{x},\sigma)$ with respect to $\boldsymbol{x}$. 
Besides, due to the convexity of $\Psi$, its Legendre conjugate $\Psi^*$ is also convex.  Therefore, the total complementary function $\Gamma(\boldsymbol{x},\sigma)$ is  concave in $\sigma$.

Therefore,  we can obtain the  global optimality of $(\boldsymbol{x}^{\Diamond}, \boldsymbol{\sigma}^{\Diamond})$ on $\boldsymbol{\Omega}\times \mathbb{E}^{+} $, i.e.,
\begin{equation*}
	\Gamma(\boldsymbol{x}^{\Diamond},  \sigma)\leq\Gamma(\boldsymbol{x}^{\Diamond},  \sigma^{\Diamond}) \leq \Gamma(\boldsymbol{x},  \sigma^{\Diamond}),\quad \forall \boldsymbol{x}\in \boldsymbol{\Omega},\; \sigma\in \mathbb{E}^{+},
\end{equation*}
which implies
$
	J_{i}(x_{i}^{\Diamond}, \boldsymbol{x}_{-i}^{\Diamond})\leq J_{i}(x_{i}, \boldsymbol{x}_{-i}^{\Diamond}), \quad \forall x_{i}\in \Omega_i,\quad \forall i\in\mathcal{N}_s.
$
Thus, $\boldsymbol{x}^{\Diamond}$ is the NE (global NE)	 of (\ref{f1}).

Next, we prove the necessity. 
Under Assumption 1,	if
	a profile $\boldsymbol{x}^{\Diamond}$ is the  NE 	of the non-convex  game (\ref{f1}), then it satisfies the first-order condition \eqref{chain} and \eqref{dde}.
 Also,
 we have  $2(\|x_{i}^{\Diamond}-x_{j}^{\Diamond}\|^2-d_{i j}^2)=0,  \forall (i, j) \in \mathcal{E}_{ss} $ and $
		2(\|x_{i}^{\Diamond}-a_{l}\|^2-e_{i l}^2)=0, \forall (i, l) \in \mathcal{E}_{as}.$
By taking $\sigma_{i j}^{s\Diamond}=2(\|x_{i}^{\Diamond}-x_{j}^{\Diamond}\|^2-d_{i j}^2)$ and $\sigma_{i l}^{a\Diamond}=2(\|x_{i}^{\Diamond}-a_{l}\|^2-e_{i l}^2)$, we obtain ${\sigma}^{\Diamond}=\boldsymbol{0}_{q}\in \mathbb{E}^{+}$. Also, the following duality relation holds:
\begin{align*}
	\sigma^{\Diamond}=\nabla \Psi(
\xi)|_{\xi=\Lambda(\boldsymbol{x}^{\Diamond})}\Longleftrightarrow \;\Lambda(\boldsymbol{x}^{\Diamond})=\nabla \Psi^{*}(\sigma^{\Diamond}). 
\end{align*}

 Thus, by substituting  $\nabla \Psi(
\Lambda(\boldsymbol{x}^{\Diamond}))$ with $\sigma^{\Diamond}$ in \eqref{chain}, we can derive that $(\boldsymbol{x}^{\Diamond},{\sigma}^{\Diamond}) \in \boldsymbol{\Omega}\times\mathbb{E}^{+} $ 
	is a stationary point of    $ \Gamma (\boldsymbol{x},\sigma)$. This completes the proof. \hfill $\square$

 \vspace{-0.13cm}

 \section{Proof of Corollary \ref{c1}}\label{cc}
For $i\in\mathcal{N}_s$, denote the unit square constraints $[x_{i,min},x_{i,max}]^{n}$ as a circumscribed polyhedron for $\Omega_i$, where $[x_{i,min},x_{i,max}]$ is the maximum diameter of $\Omega_i$.
Then we can describe the bound  $W$ as
$$
\begin{aligned}
	W=\operatorname{min}\{&{\operatorname{min}}_{(i,j)\in\mathcal{E}_{ss}}\{2\|x_{i,max}-x_{i,min}\|^2-2d_{ij}^2\},\\
	&  {\operatorname{min}}_{(i,l)\in\mathcal{E}_{as}}\{2\|x_{i,max}-a_{l}\|^2-2e_{il}^2\}\}.
\end{aligned}
$$ 
Accordingly, the 
feasible range of $\sigma_{i j}^s$ is 
converted into $[0,W]$, as well as 
$\sigma_{i l}^a$.
Recall that the  global NE
can be represented as  a strategy profile  $\boldsymbol{x}^{\Diamond}$ that satisfies $\|x_{i}^{\Diamond}-x_{j}^{\Diamond}\|^2-d_{i j}^2=0$ and $\|x_{i}^{\Diamond}-a_{l}\|^2-e_{i l}^2=0$ for any $(i,j)\in \mathcal{E}_{ss}$ and $(i,l)\in \mathcal{E}_{as}$.
It can be deduced from  the dual relation (\ref{dual}) that $\sigma_{i j}^{s\Diamond}=0$ and $\sigma_{i l}^{a\Diamond}=0$,
as shown in \eqref{zero}.
Therefore, $\sigma^{\Diamond}=\operatorname{col}\{\operatorname{col}\{\sigma_{ij}^{s\Diamond}\}_{(i,j)\in\mathcal{E}_{ss} },\operatorname{col}\{\sigma_{il}^{a\Diamond}\}_{(i,l)\in\mathcal{E}_{as} }
\} \in [0,W]^q$.
Similar to the proof in Theorem \ref{t2}, we get the result. \hfill $\square$

\vspace{-0.3cm}


\section{Proof of Theorem  \ref{cent}}\label{bbw}
According to \eqref{cen-compact},
\begin{align}\label{zeta}
	\|\boldsymbol{\zeta}[k+1]-\boldsymbol{\zeta}^{\Diamond}
	&=\| \Pi_{\boldsymbol{\Xi}}\left(\boldsymbol{\zeta}[k]-\alpha[k]F(\boldsymbol{\zeta}[k])\right)-\boldsymbol{\zeta}^{\Diamond} \|_{2}^{2}\nonumber\\
	&\leq \| \boldsymbol{\zeta}[k]-\boldsymbol{\zeta}^{\Diamond}-\alpha[k]F(\boldsymbol{\zeta}[k])\|_{2}^{2}\nonumber\\
&\leq \| \boldsymbol{\zeta}[k]\!-\!\boldsymbol{\zeta}^{\Diamond} \|_{2}^{2}\!-\!2\alpha[k](\boldsymbol{\zeta}[k]\!-\!\boldsymbol{\zeta}^{\Diamond})^{\operatorname{T}}F(\boldsymbol{\zeta}[k])\!\nonumber\\
	&\quad+ \!\alpha[k]^{2}\| F(\boldsymbol{\zeta}[k])\|_{2}^{2}. 
\end{align}
Moreover,  due to $\sigma\!\in\mathbb{E}^{+}$ in (\ref{s23}) and the   convexity of $\Psi^{*}$, 
	\begin{align}\label{rtyt}
		&\quad\langle (\boldsymbol{\zeta}[k]), \boldsymbol{\zeta}^{\Diamond}-\boldsymbol{\zeta}[k]\rangle\notag \\
		&\leq \sigma[k]^{\mathrm{T}} \Lambda(\boldsymbol{x}^{\Diamond})-\Psi^{*}\sigma[k])-(\sigma^{\Diamond \mathrm{T}} \Lambda(\boldsymbol{x}[k]^{\mathrm{T}})-\Psi^{*}(\sigma^{\Diamond}))\notag\\
		&= \Gamma(\boldsymbol{x}^{\Diamond},\sigma[k])-\Gamma(\boldsymbol{x}[k],\sigma^{\Diamond}).\vspace{-0.2cm}
	\end{align}
 By substituting  \eqref{rtyt} into \eqref{zeta} and rearranging the terms therein, we have 
\begin{align*}
	&\alpha[k]\!(\Gamma(\boldsymbol{x}[k]\!,\sigma^{\Diamond}\!)\!-\!\Gamma(\boldsymbol{x}^{\Diamond}\!,\sigma[k]))\!\leq\!\alpha[k]\! F(\boldsymbol{\zeta}[k])^{\mathrm{T}}\!(\boldsymbol{\zeta}[k]\!-\!\boldsymbol{\zeta}^{\Diamond}\!)\\
	&\quad\quad\leq \frac{1}{2}\| \boldsymbol{\zeta}[k]-\boldsymbol{\zeta}^{\Diamond} \|_{2}^{2}\!-\!\frac{1}{2}\| \boldsymbol{\zeta}[k+1]\!-\!\boldsymbol{\zeta}^{\Diamond} \|_{2}^{2}\!+ \!\frac{\alpha[k]^2}{2}\|F(\boldsymbol{\zeta}[k]) \|^2.
\end{align*}
Consider the sum of the above inequalities over $1, \dots, k$,
\begin{align}\label{o1}
	&\sum_{j=1}^k \alpha[j](
	\Gamma(\boldsymbol{x}[j],\sigma^{\Diamond})\!-\!
	\Gamma(\boldsymbol{x}^{\Diamond}\!,\sigma[j])\!
	) \\
 &\quad\quad\quad\leq\frac{ \| \boldsymbol{\zeta}[1]-\boldsymbol{\zeta}^{\Diamond} \|_{2}^{2}\!+\!\sum^{k}_{j=1}\!\!\alpha[j]^2 \!{M}^2_{2}}{2},\notag\vspace{-0.6cm}
\end{align}
where $\|F(\boldsymbol{\zeta})\|^2$ is bounded by a constant $M_2$.
	Denote the weight by $\lambda[j]= {\alpha[j]}/{\sum\nolimits^{k}_{l=1}\alpha[l]}$. Then (\ref{o1})  yields
	\begin{align*}
		&\quad\,{(\frac{1}{2}\| \boldsymbol{\zeta}[1]-\boldsymbol{\zeta}^{\Diamond} \|_{2}^{2} +\frac{1}{2}{M}^2_{2}\sum\nolimits^{k}_{j=1}\alpha_{j}^2 )}/{\sum\nolimits^{k}_{j=1}\alpha[j]}\\	
		&\geq\sum\nolimits_{j=1}^{k} \frac{\alpha[j]}{\sum^{k}_{l=1}\alpha_{l}} \left(
		\Gamma(\boldsymbol{x}[j],\sigma^{\Diamond})-
		\Gamma(\boldsymbol{x}^{\Diamond},\sigma[j])
		 \right)\\  
		&= \sum\nolimits_{j=1}^{k}  \lambda[j] \left(
		\Gamma(\boldsymbol{x}[j],\sigma^{\Diamond})-
		\Gamma(\boldsymbol{x}^{\Diamond},\sigma[j])
		\right)\\
		&\geq \Gamma(\sum\nolimits_{j=1}^{k}  \lambda[j] \boldsymbol{x}[j],\sigma^{\Diamond})-
		\Gamma(\boldsymbol{x}^{\Diamond},\sum\nolimits_{j=1}^{k}  \lambda[j] \boldsymbol{\sigma}[j])\\
		&= \Gamma(\hat{\boldsymbol{x}}[k],\sigma^{\Diamond})-
		\Gamma(\boldsymbol{x}^{\Diamond},\hat{\sigma}[k]),
	\end{align*}
	where the last inequality is true due to Jensen's inequality. 
	Since the step size satisfies $\alpha[k]=\sqrt{2{d}}/{M}_2\sqrt{k}$, we finally derive that 
	\begin{align*}
		\Gamma(\hat{\boldsymbol{x}}[k],\sigma^{\Diamond})-
		\Gamma(\boldsymbol{x}^{\Diamond},\hat{\sigma}[k])
		\leq \sqrt{{2d}}{M}_2/{\sqrt{k}} .
  \vspace{-0.4cm}
	\end{align*}
	Moreover, if $	\sigma_{i j}^{s\Diamond}=2(\|x_{i}^{\Diamond}-x_{j}^{\Diamond}\|^2-d_{i j}^2)$, and 
	$\sigma_{i l}^{a\Diamond}=2(\|x_{i}^{\Diamond}-a_{l}\|^2-e_{i l}^2)$ for any $  (i, j) \in \mathcal{E}_{ss}$ and $(i, l) \in \mathcal{E}_{as}$, then the convergent point $\boldsymbol{x}^{\Diamond}$ is indeed the global NE. 
\hfill $\square$
\vspace{-0.2cm}
\section{Detailed  design ideas for  Algorithm 2}\label{ccf}

We show more details about the design $({\sigma_{ij}^{s_{i}}+\sigma_{ij}^{s_{j}}})/{2}$ in Algorithm 2. As mentioned before,   the design is inspired by sliding mode control \cite{sarpturk1987stability}, where  $\mathcal{J}_{\boldsymbol{F}}({\boldsymbol{z}[k]}) $ is viewed as the sliding manifold.  
 Specifically,  if   we follow the  existing distributed methods  by only employing $\sigma_{ij}^{s_i}$ to decompose $\Gamma(\boldsymbol{x},\sigma)$, 
 the  local complementary problem for non-anchor node $i$ is described as
 \begin{equation*}
 	\begin{aligned}
 		&\quad\;{\varGamma}_{i}^{\prime}(x_{i},\{x_j\}_{j\in\mathcal{N}_{s}^{i} },\sigma_i)\\
 		&=\!\sum\nolimits_{j\in\mathcal{N}_{s}^{i}  }\!\sigma_{ij}^{s_{i}}(\|x_{i}\!-\!\!x_{j}\|^{2}-\!d_{ij}^2)\!+\!\sum\nolimits_{l\in\mathcal{N}_{a}^{i}  }\!\sigma_{ij}^{a_{i}}(\|x_{i}\!-\!a_{l}\|^{2}\!-\!e_{il}^2)\\
 		&\quad-\sum\nolimits_{j\in\mathcal{N}_{s}^{i}  }{(\sigma_{ij}^{s_{i}})^2}/{4}\!-\!\sum\nolimits_{l\in\mathcal{N}_{a}^{i}}{(\sigma_{il}^{a_{i}})^2}/{4}. 
 	\end{aligned}
 \end{equation*}
 Denote the pseudo-gradient of ${\varGamma}_{i}^{\prime}$ for $i\in \mathcal{N}_s$ by  $${\boldsymbol{F}}^{\prime}(\boldsymbol{z})=\operatorname{col}\{\operatorname{col}\{\nabla_{x_{i}}{\varGamma}^{\prime}_{i}\}_{i=1}^N, \operatorname{col}\{-\nabla_{\sigma_{i}}{\varGamma}^{\prime}_{i}\}_{i=1}^N\},$$ where 
 $\nabla_{x_{i}}{\varGamma}^{\prime}_{i}=\!\!\!\!\sum_{j\in\mathcal{N}_{s}^{i} }\!2\sigma_{ij}^{s_{i}}(x_{i}-x_{j})\!\!+\!\!\sum_{l\in\mathcal{N}_{a}^{i} }\!\!2\sigma_{il}^{a_{i}}\!(x_{i}\!-\!a_{l})$ and 
 $$\nabla_{\sigma_{i}}{\varGamma}^{\prime}_{i}=
 \left[\begin{array}{c}
 	\operatorname{col}\{\|x_{i}-x_{j}\|^{2}\!-\!d_{ij}^2\!-\!\frac{1}{2}\sigma_{ij}^{s_{i}} \}_{j\in\mathcal{N}_{s}^{i} }\\ \operatorname{col}\{\|x_{i}\!-\!a_{l}\|^{2}-e_{il}^2\!-\!\frac{1}{2}\sigma_{il}^{a_{i}}\}_{l\in\mathcal{N}_{a}^{i} }
 \end{array}\right].$$
Denote the  Jacobian matrix of  $\boldsymbol{F}({\boldsymbol{z}[k]}) $ by
 $\mathcal{J}_{{\boldsymbol{F}}^{\prime}}({\boldsymbol{z}[k]}) $, where 
  $$
 \mathcal{J}_{{\boldsymbol{F}}^{\prime}}({\boldsymbol{z}})=\left[\begin{array}{cc}
 	A &B  \\
 	C & D
 \end{array}\right],\vspace{-0.2cm}
 $$
 with
 $$
 A\!\!=\!\!\left[\begin{array}{ccc}
 	\!\!\!\!\nabla_{x_1 x_1}^2\!\! {\varGamma}^{\prime}_1\!\! \!\!&\!\!\! \!\dots\!\!\! \!\!& \!\nabla_{x_1 x_N}^2\! {\varGamma}^{\prime}_1 \!\!\!\! \\
\!\! \!\!	\vdots\!\! \!\!&\!\! \!\!\ddots\!\! &\!\! \!\!\vdots\!\!\!\!\\
 \!\!\!\!	\nabla_{x_N x_1}^2\!\! {\varGamma}^{\prime}_N \!\!\!\!&\!\! \!\!\dots\!\!& \!\!\!\nabla_{x_N x_N}^2 \!\! {\varGamma}^{\prime}_N\!\!\!\!
 \end{array}\right]\!, B\!\!=\!\!\left[\begin{array}{ccc}
 \!\!	\!\!\nabla_{ x_1\sigma_1}^2\!\! {\varGamma}^{\prime}_1\!\!\!\! & \!\!\!\! \dots\!\! \!\!& \nabla_{x_1\sigma_N}^2\!\! {\varGamma}^{\prime}_1 \!\!\!\!\\
 \!\! \!\!	\vdots \!\!\!\!&\!\!\!\! \ddots \!\!\!\!&\!\!\!\! \vdots\!\!\!\!\\
 \!\! \!\!	\nabla_{x_N \sigma_1}^2\!\! {\varGamma}^{\prime}_N\!\!\!\! &\!\! \!\!\dots \!\!\!\!&\!\!\!\! \nabla_{x_N \sigma_N}^2\!\! {\varGamma}^{\prime}_N \!\!\!\! \end{array}\right]\!,
 $$
 $$
 C\!\!=\!\!\left[\begin{array}{ccc}
  \!\!\!\! \!\!	-\!\nabla_{\sigma_1 x_1}^2 \!\! {\varGamma}^{\prime}_1\!\!\!\! &\!\!\!\! \dots\!\!\!\!&\!\!\!\! -\nabla_{\sigma_1 x_N}^2  {\varGamma}^{\prime}_1    \!\! \!\!\!\\
 \!\! \!\! \!\!	\vdots\!\!\!\! & \!\!\!\!\ddots\!\!\!\! &\!\! \!\!\vdots \!\! \!\!\!\\
 	\!\! \!\! \!-\!\nabla_{\sigma_N x_1}^2 \!\! {\varGamma}^{\prime}_N \!\!\!\!&\!\! \!\!\dots\!\!\!\!&\!\!\! -\!\nabla_{\sigma_N x_N}^2 \!\! {\varGamma}^{\prime}_N  \!\! \!\!\!
 \end{array}\right]\!, D\!\!=\!\!\left[\begin{array}{ccc}
 \!\!\!\!	-\!\nabla_{ \sigma_1 \sigma_1}^2 \!\!{\varGamma}^{\prime}_1 \!\!\!\!&\!\! \!\!\dots\!\!\!\! &\!\!\!\!\! -\!\nabla_{\sigma_1 \sigma_N }^2 \!\!{\varGamma}^{\prime}_1  \!\!\!\! \!\!\\
 \!\!\!\!	\vdots \!\!\!\!&\!\!\!\! \ddots \!\!\!\!&\!\!\!\! \vdots\!\!\!\!\!\!\\
 \!\!\!\!	-\!\nabla_{\sigma_N \sigma_1}^2\!\! {\varGamma}^{\prime}_N \!\!\!\!\!&\!\! \dots\!\!\!\! & \!\!\!-\!\nabla_{\sigma_N \sigma_N}^2 \!\!{\varGamma}^{\prime}_N  \!\!\!\!\!\!
 \end{array}\right]\!.
 $$
 The form  $\mathcal{J}_{\boldsymbol{F}}({\boldsymbol{z}[k]}) $  can be expressed in a simliarly way. However, we cannot guarantee the symmetricity and positive semidefiniteness of $\frac{1}{2}(\mathcal{J}_{{\boldsymbol{F}}^{\prime}}({\boldsymbol{z}[k]})^T+\mathcal{J}_{{\boldsymbol{F}}^{\prime}}({\boldsymbol{z}[k]}) )$, since we only have $B^T\neq-C$ in such a distributed case. 
 Thus, our goal is to control   $ {\boldsymbol{F}^{\prime}}(\boldsymbol{z}[k])+ \boldsymbol{u}[k]$ to ensure that $\frac{1}{2}(\mathcal{J}_{{\boldsymbol{F}}^{\prime}}({\boldsymbol{z}[k]})^T+\mathcal{J}_{{\boldsymbol{F}}^{\prime}}({\boldsymbol{z}[k]})) $ lies with a space of  positive semidefinite matrices.
  Accordingly,  by substituting $\sigma_{ij}^{s_{i}}$ with the average  $({\sigma_{ij}^{s_{i}}+\sigma_{ij}^{s_{j}}})/{2}$ and modifying ${\varGamma}^{\prime}_i$ to $\varGamma_{i}$, we verify that $\mathcal{J}_{\boldsymbol{F}}$ on ${\boldsymbol{\varXi}}$ satisfies 
 $$
 \frac{1}{2}(\mathcal{J}_{\boldsymbol{F}}^T+\mathcal{J}_{\boldsymbol{F}})=\left[\begin{array}{cc}
 	A &0  \\
 	0 & D
 \end{array}\right] \succeq 0.\vspace{-0.2cm}
 $$
 In this light, 
 $\mathcal{J}_{\boldsymbol{F}}({\boldsymbol{z}[k]}) $  is regarded as the sliding manifold and
  the  control $\boldsymbol{u}[k]$ is constructed as in \eqref{contro}.
Here, we select the initial value which is located on the sliding manifold for convenience, i.e.,   taking $\sigma_{ij}^{s_{i}}[0]=\sigma_{ij}^{s_{j}}[0]$.

\vspace{-0.2cm}

\section{Proof of Theorem \ref{t4}}\label{c2f}
Suppose that $(\tilde{\boldsymbol{z}}, {\boldsymbol{z}})$ is a fixed point of  (\ref{dis-compact}), that is
\begin{equation}\label{fix}
	\left\{\begin{array}{l}
		\tilde{\boldsymbol{z}}=\Pi_{{\boldsymbol{\varXi}}}(\boldsymbol{z}-\beta \boldsymbol{F}(\boldsymbol{z})), \\
		\boldsymbol{z}=\Pi_{{\boldsymbol{\varXi}}}(\boldsymbol{z}-\beta \boldsymbol{F}(\tilde{\boldsymbol{z}})).
	\end{array}\right.\vspace{-0.2cm}
\end{equation}
Then 
\begin{equation*}
	\begin{aligned}
		\|\tilde{\boldsymbol{z}}-\boldsymbol{z}\| & =\left\|\Pi_{\boldsymbol{\varXi}}[\boldsymbol{z}-\beta \boldsymbol{F}(\boldsymbol{z})]-\Pi_{\boldsymbol{\varXi}}[\boldsymbol{z}-\beta \boldsymbol{F}(\tilde{\boldsymbol{z}})]\right\| \\
		& \leq \beta\|\boldsymbol{F}(\boldsymbol{z})-\boldsymbol{F}(\tilde{\boldsymbol{z}})\|\leq\beta L\|\tilde{\boldsymbol{z}}-\boldsymbol{z}\| 
	\end{aligned}
\end{equation*}
which indicates
$
(1-\beta\cdot L)\|\tilde{\boldsymbol{z}}-\boldsymbol{z}\|\leq0
$. Since the stepsize satisfies  $\beta\cdot L<1$, 
\begin{equation}\label{fix2}
	\tilde{\boldsymbol{z}}=\boldsymbol{z}=\Pi_{{\boldsymbol{\varXi}}}(\boldsymbol{z}-\beta \boldsymbol{F}(\boldsymbol{z})).
\end{equation}
Take $\boldsymbol{z}=\operatorname{col}\{\boldsymbol{x}^{\Diamond},\boldsymbol{\sigma}^{\Diamond}\}$.  Then (\ref{fix2}) is rewritten as 
\begin{subequations}\label{e346578}
	\begin{align}
		&(\sum\nolimits_{j\in\mathcal{N}_{s}^{i} }\!(\sigma_{ij}^{s_{i}\Diamond}+\sigma_{ij}^{s_{j}\Diamond})(x_{i}^{\Diamond}\!-\!x_{j}^{\Diamond})\!\!\\
  &\quad\quad+\!\!\!\sum\nolimits_{l\in\mathcal{N}_{a}^{i} }\!\!\!2\sigma_{il}^{a_{i}}\!(x_{i}^{\Diamond}\!-\!a_{l}))^{\operatorname{T}} \!(x_{i}\!-\!x_{i}^{\Diamond})\!\geq\!0, \quad \forall x_{i}\in \Omega_i,\notag\\
		&(-(\|x_{i}^{\Diamond}-x_{j}^{\Diamond}\|^{2}\!-\!d_{ij}^2)/2\!+\!({\sigma_{ij}^{s_{i}\Diamond}\!\!+\sigma_{ij}^{s_{j}\Diamond}})/{8})^{\operatorname{T}} \notag\\
  &\quad\quad\quad\quad\cdot(\sigma_{ij}^{s_{i}}-\sigma_{ij}^{s_{i}\Diamond})\geq0, \quad\forall \sigma_{ij}^{s_{i}}\in \mathscr{E}_{i}^+, \,j\in\mathcal{N}_{s}^{i},\\
		& (-(\|x_{i}^{\Diamond}-a_{l}\|^{2}\!-\!e_{ij}^2)\!+\!{\sigma_{il}^{a_{i}\Diamond}}/{2})^{\operatorname{T}}\notag\\
  &\quad\quad\quad\quad\cdot(\sigma_{il}^{a_{i}}-\sigma_{il}^{a_{i}\Diamond})\geq0, \quad\forall \sigma_{il}^{a_{i}}\in \mathscr{E}_{i}^+, \,l\in\mathcal{N}_{a}^{i}.
	\end{align}
\end{subequations}
Since $\sigma_{ij}^{s_{j}\Diamond}=\sigma_{ij}^{s_{i}\Diamond}$, (\ref{e346578})(a)-(b) becomes  
\begin{subequations}\label{e346578p}
\begin{align}
		&(\sum\nolimits_{j\in\mathcal{N}_{s}^{i} }\!2\sigma_{ij}^{s_{i}\Diamond}(x_{i}^{\Diamond}-x_{j}^{\Diamond})\\
  &\quad\quad+\sum\nolimits_{l\in\mathcal{N}_{a}^{i} }\!\!\!2\sigma_{il}^{a_{i}}\!(x_{i}^{\Diamond}\!-\!a_{l}))^{\operatorname{T}} (x_{i}-x_{i}^{\Diamond})\geq0, \quad \forall x_{i}\in \Omega_i,\notag\\
	&(-(\|x_{i}^{\Diamond}-x_{j}^{\Diamond}\|^{2}\!-\!d_{ij}^2)/2\!+\!{\sigma_{ij}^{s_{i}\Diamond}\!\!}/{4})^{\operatorname{T}}\\
  &\quad\quad\quad\quad\cdot(\sigma_{ij}^{s_{i}}-\sigma_{ij}^{s_{i}\Diamond})\geq0, \quad\forall \sigma_{ij}^{s_{i}}\in \mathscr{E}_{i}^+, \,j\in\mathcal{N}_{s}^{i}.\notag \vspace{-0.4cm}
\end{align}
\end{subequations}
{Since $\varGamma_i$ is strictly concave in $\sigma_i$, $\sigma_i^{\Diamond}=\boldsymbol{0}_{q_i}$. Moreover, it follows from the duality relation that $\sigma_{i j}^{s_{i}\Diamond}\!=\!2(\|x_{i}^{\Diamond}-x_{j}^{\Diamond}\|^2-d_{i j})\!$ and $
	\sigma_{i l}^{a_{i}\Diamond}\!=\!2(\|x_{i}^{\Diamond}-a_{l}\|^2-e_{i l}^2)$.
}
Thus, according to the chain rule, (\ref{e346578p})(a) becomes
$
	\mathbf{0}_{n} \in \nabla_{x_{i}} J_{i}(x_{i}^{\Diamond}, \boldsymbol{x}_{-i}^{\Diamond})+\mathcal{N}_{\Omega_{i}}({x_{i}^{\Diamond}}). 
$
Moreover, when  $ \sigma_{i}\in \mathscr{E}_{i}^{+}$,  the Hessian matrix satisfies
\begin{align*}
	\nabla_{x_i}^{2} \varGamma_{i}
	=\sum\nolimits_{j\in\mathcal{N}_{s}^{i}  }(\sigma_{ij}^{s_{i}}\!+\!\sigma_{ij}^{s_{j}})\boldsymbol{I}_{n}+\sum\nolimits_{l\in\mathcal{N}_{a}^{i}  }\!\!2\sigma_{iL}^{a_{i}} \boldsymbol{I}_{n}\succeq \boldsymbol{0}_n,
\end{align*} 
which indicates 
the convexity of $\varGamma_{i}(x_{i},\{x_j\}_{j\in\mathcal{N}_{s}^{i} },\sigma_i,$ $\!\{\sigma_j\}_{j\in\mathcal{N}_{s}^{i} })$ with respect to $x_{i}$. 
Besides, $\varGamma_{i}(x_{i},\!\{x_j\}_{j\in\mathcal{N}_{s}^{i} },$ $\sigma_i,\!\{\sigma_j\}_{j\in\mathcal{N}_{s}^{i} })$ is concave in  $\sigma_{i}$.  
In this light, we can obtain the global optimality of $(\boldsymbol{x}^{\Diamond}, \boldsymbol{\sigma}^{\Diamond})$ on $\boldsymbol{\Omega}\times \bar{\mathscr{E}}^{+} $, i.e.,
for $i\in\mathcal{N}_s$,
\begin{equation*}
	\begin{aligned}
		&\varGamma_{i}(x_{i}^{\Diamond},\!\{\!x_j^{\Diamond}\!\}_{j\in\mathcal{N}_{s}^{i} },\sigma_i,\!\{\!\sigma_j^{\Diamond}\!\}_{j\in\mathcal{N}_{s}^{i} })\!\!\leq\!\!\varGamma_{i}(x_{i}^{\Diamond},\!\{\!x_j^{\Diamond}\!\}_{j\in\mathcal{N}_{s}^{i} },\sigma_i^{\Diamond},\!\{\!\sigma_j^{\Diamond}\!\}_{j\in\mathcal{N}_{s}^{i} }) \\
		&\leq \varGamma_{i}(x_{i},\!\{\!x_j^{\Diamond}\!\}_{j\in\mathcal{N}_{s}^{i} },\sigma_i^{\Diamond},\!\{\!\sigma_j^{\Diamond}\!\}_{j\in\mathcal{N}_{s}^{i} }),\quad \forall x_{i}\in \Omega_i,\; \sigma_{i}\in \mathscr{E}_{i}^{+}, 	
	\end{aligned}
\end{equation*}
which implies
$
	J_{i}(x_{i}^{\Diamond}, \boldsymbol{x}_{-i}^{\Diamond})\leq J_{i}(x_{i}, \boldsymbol{x}_{-i}^{\Diamond}), \quad \forall x_{i}\in \Omega_i,\quad \forall i\in\mathcal{N}_s.
$
Thus, $\boldsymbol{x}^{\Diamond}$ is the global NE of (\ref{f1}).

Conversely, suppose that  $\boldsymbol{x}^{\Diamond}$ is the global NE of (\ref{f1}). By Theorem \ref{t2} and Corollary \ref{c1}, there exists $\sigma^\Diamond=\operatorname{col}\{\operatorname{col}\{\sigma_{ij}^{s\Diamond}\}_{(i,j)\in\mathcal{E}_{ss}},\operatorname{col}\{\sigma_{il}^{a\Diamond}\}_{(i,l)\in\mathcal{E}_{as}} \}\in\mathbb{R}^{q}$  such that   the duality relation \eqref{dual} is satisfied and 
\begin{equation*}
	\begin{aligned}
 &(\sum\nolimits_{j\in\mathcal{N}_{s}^{i} }\!2\sigma_{ij}^{s\Diamond}(x_{i}^{\Diamond}\!-\!x_{j}^{\Diamond})\!\!\\\vspace{-0.2cm}
  &\quad\quad+\!\!\!\sum\nolimits_{l\in\mathcal{N}_{a}^{i} }\!\!\!2\sigma_{il}^{a}\!(x_{i}^{\Diamond}\!-\!a_{l}))^{\operatorname{T}} \!(x_{i}\!-\!x_{i}^{\Diamond})\!\geq\!0, \quad \forall x_{i}\in \Omega_i,\vspace{-0.2cm}
  	\end{aligned}
\end{equation*}
\begin{equation*}
	\begin{aligned}
		&(-(\|x_{i}^{\Diamond}-x_{j}^{\Diamond}\|^{2}\!-\!d_{ij}^2)\!+\!{\sigma_{ij}^{s\Diamond}}/{2})^{\operatorname{T}} \\
  &\quad\quad\quad\cdot(\sigma_{ij}^{s}-\sigma_{ij}^{s\Diamond})\geq0, \quad\forall \sigma_{ij}^{s}\in [0,W],\,j\in\mathcal{N}_{s}^{i}, \\\vspace{-0.2cm}
		& (-(\|x_{i}^{\Diamond}-a_{l}\|^{2}\!-\!e_{ij}^2)\!+\!{\sigma_{il}^{a\Diamond}}/{2})^{\operatorname{T}}\\
  &\quad\quad\quad\cdot(\sigma_{il}^{a}-\sigma_{il}^{a\Diamond})\geq0, \quad\forall \sigma_{il}^{a}\in [0,W], \,l\in\mathcal{N}_{a}^{i}.
	\end{aligned}
\end{equation*}
For $i\in\mathcal{N}_s$, by
choosing ${\sigma}_i^{\Diamond}=\operatorname{col}\{\operatorname{col}\{\sigma_{ij}^{s_{i}\Diamond}\}_{j\in\mathcal{N}_{s}^{i} },$ $\operatorname{col}\{\sigma_{il}^{a_{i}\Diamond}\}_{l\in\mathcal{N}_{a}^{i} }\}$ as $\sigma_{ij}^{s_{i}\Diamond}=\sigma_{ij}^{s_{i}\Diamond}=\sigma_{ij}^{s\Diamond}$, $\sigma_{il}^{a_{i}\Diamond}=\sigma_{il}^{a\Diamond}$, 
one can verify that $(\boldsymbol{x}^{\Diamond},\boldsymbol{\sigma}^{\Diamond})$ renders \eqref{dual-rela} and \eqref{e346578}, which completes the proof.
\hfill $\square$

\end{appendices}
\vspace{-0.2cm}
	\bibliographystyle{IEEEtran}
\bibliography{reference}

\end{document}